\documentclass[12pt]{article} 
\author{Constantin N. Beli}
\title{Four conjectures by Zhi-Hong Sun}
\usepackage{amsfonts,amsmath,amssymb}
\def\j{\infty} 
    
\def\p{\partial} \def\q{\nabla}  \def\s{\sigma}

  \def\z{\longrightarrow}
\def\({\overline} \def\){\underline}
\def\<{\cdot} \def\go{\mathfrak} \def\>{~~~~~~~} \def\#{{\bf
Definition}} \def\*{\section} \def\be{\begin{equation}}
\def\ee{\end{equation}}

\newcommand{\legendre}[2]{\genfrac{(}{)}{}{}{#1}{#2}}

\def\sb{\subset}  \def\sbq{\subseteq} \def\la{\langle}
\def\ra{\rangle}  \def\ti{\times} 
\def\od{{\rm ord}\,}\def\oo{{\cal O}}

\def\p{\go p} \def\q{\go q} \def\P{\go P} \def\Q{\go Q} \def\mm{\go m}
\def\*{\sharp}   \def\0{} 
 \def\1{^{-1}}  
\def\bmat{\left(\begin{array}} \def\emat{\end{array}\right)}
\def\ev{\equiv} \def\ap{\cong}  
 \def\N{{\rm N}} \def\CC{{\cal C}}
\def\no{\noindent}
 \def\btm{\begin{thm}}
\def\etm{\end{tm}}
 \def\blem{\begin{lem}}
\def\elem{\end{lem}}

\newtheorem{theorem}{Theorem}[section]

\newtheorem{lemma}[theorem]{Lemma}
\newtheorem{definition}{Definition}
\newtheorem{corollary}[theorem]{Corollary}

\newtheorem{bof}[theorem]{}
\newtheorem{teorema}{Theorem}

\def\qed{\mbox{$\Box$}\vspace{\baselineskip}}
\def\pf{$Proof.$} 
\def\bco{\begin{corollary}} \def\eco{\end{corollary}} 
\def\bdf{\begin{definition}} \def\edf{\end{definition}} 
\def\btm{\begin{theorem}} \def\etm{\end{theorem}} 
\def\blm{\begin{lemma}} \def\elm{\end{lemma}} 
\def\bff{\begin{bof}\rm} \def\eff{\end{bof}}
\def\btr{\begin{teorema}} \def\etr{\end{teorema}}

\def\de{\newcommand} \de\tm[1]{{\no\bf Theorem~#1}} 
 \def\mb{\mathbb} 
\def\RR{{\mb R}}\def\QQ{{\mathbb Q}}\def\CC{{\mb C}}\def\ZZ{{\mb
Z}} 

\begin{document}
\maketitle

\begin{quote}
{\bf\footnotesize We prove some results conjectured by Zhi-Hong Sun
regarding the value of $\varepsilon_d^{\frac{p-1}4}\mod p$, where $\varepsilon_d$ is a unit
of norm $-1$ in $\QQ (\sqrt d)$, $\legendre dp=1$ and $p\ev 1\pmod 4$. The
answer is given in terms of how we write $p$ as $p=f(x,y)=u^2+v^2$,
with $x,y,u,v\in\ZZ$,  where $f$  is a certain quadratic form of
determinant $-4d$. In Sun's conjectures $d$ has the particular form
$d=b^2+4$ with $b$ odd or $d=b^2/4+1$ with $b$ even. Also $f$ has a
particular form, namely $f=X^2+dY^2$. We will show how this problem
can be tackled for $d>1$ squarefree arbitrary and $f$ arbitrary with
discriminant $-d$ or $-4d$.}

\end{quote}


\section{Introduction}

Let $d>1$ be a square-free integer and let $\varepsilon$ be an integer in $\QQ
(\sqrt d)$. In [S1] Zhi-Hong Sun determined the value of
$\varepsilon^{(p+\legendre{-1}p)/2}\mod p$ in the case when $\legendre{-d}p=1$. His
answer is given in terms of $x,y\in\ZZ$ satisfying $p=f(x,y)$, where
$f=AX^2+BXY+CY^2$ is a quadratic form of discriminant $B^2-4AC=-k^2d$
and $k$ is a bounded positive integer. Such a representation by a
quadratic form exists because $\legendre{-d}p=1$. Later Sun stated several
conjectures regarding the value of $\varepsilon^{(p+\legendre{-1}p)/4}\mod p$, again
when $\legendre{-d}p=1$. In [B] we proved two of his conjectures in the
case when $\legendre dp=\legendre{-1}p=-1$. The value of $\varepsilon^{(p+1)/4}\mod p$
can be determined in terms of $x,y,u,v\in\ZZ$ satisfying
$p=f(x,y)=u^2\pm 2v^2$, where $f$ is a quadratic form of discriminant
$-k^2d$ and the $\pm$ sign is $+$ if $p\ev 3\pmod 8$ and it is $-$ if
$p\ev 7\pmod 8$. In some cases, such as in Sun's conjectures,
$\varepsilon^{(p+1)/4}\mod p$ can be determined in terms of $x,y$ alone but at
this time we don't know if this can be achieved always. 

The conjectures we prove in this paper deal with the value of
$\varepsilon^{(p-1)/4}\mod p$ in the case when $\legendre dp=\legendre{-1}p=1$ and the
answer is given in terms of $x,y,u,v\in\ZZ$ satisfying
$p=f(x,y)=u^2+v^2$, where $f$ is a quadratic form of discriminant
$-k^2d$. 

We define $U_n(b,c)$ and $V_n(b,c)$ to be the Lucas sequences given
by: 

$U_0(b,c)=0,~U_1(b,c)=1\text{ and
}U_{n+1}(b,c)=bU_n(b,c)-cU_{n-1}(b,c)$

\no and

$V_0(b,c)=2,~V_1(b,c)=b\text{ and
}V_{n+1}(b,c)=bV_n(b,c)-cV_{n-1}(a,b).$

Alternatively, if $d:=b^2-4c$ is not a square, the integers
$U_n(b.c),V_n(b,c)$ can be defined by the relation
$\frac{V_n(b,c)+U_n(b,c)\sqrt d}2=\left(\frac{b+\sqrt d}2\right)^n$. 

Note that $(-b)^2-4c=b^2-4c=d$ and so $\frac{V_n(-b,c)+U_n(-b,c)\sqrt
d}2=\left(\frac{-b+\sqrt d}2\right)^n=(-1)^n\left(\frac{b-\sqrt
d}2\right)^n=(-1)^n\frac{V_n(b,c)-U_n(b,c)\sqrt d}2$. Hence
$V_n(-b,c)=(-1)^nV_n(b,c)$ and $U_n(-b,c)=(-1)^{n+1}U_n(b,c)$. 

We now state the Conjectures 9.4, 9.11, 9.14 and 9.17 from [S2], which
we will call Conjectures 1-4.\footnote{In order that our notation
matches that from [B], we took the liberty of replacing Sun's $c$ and
$d$ by $u$ and $v$. Thus $p=c^2+d^2$ from [S2] becomes $p=u^2+v^2$.}

{\bf Conjecture 1} Let $p\ev 1\pmod 4$ be a prime, $b\in\ZZ$, $2\nmid
b$ and $p=x^2+(b^2+4)y^2=u^2+v^2\neq b^2+4$ with
$x,y,u,v\in\ZZ$. Suppose that $u$ and all odd parts of $x,y,v$ are
$\ev 1\pmod 4$.

(i) If $4\nmid xy$, then $V_{\frac{p-1}4}(b,-1)\ev 0\pmod p$ and 
$$U_{\frac{p-1}4}(b,-1)\ev\begin{cases}(-1)^{\frac v4}\frac{2y}x\pmod
p&\text{if }2\| x\text{ and }b\ev 1,3\pmod 8,\\ 
-(-1)^{\frac v4}\frac{2y}x\pmod p&\text{if }2\| x\text{ and }b\ev 5,7\pmod
8,\\ \frac{2yv}{xu}\pmod p&\text{if }2\| y.\end{cases}$$

(ii) If $4\mid xy$, then $U_{\frac{p-1}4}(b,-1)\ev 0\pmod p$ and
$$V_{\frac{p-1}4}(b,-1)\ev\begin{cases}2(-1)^{\frac{y+v}4}\pmod p&\text{if
}4\mid y,\\ -2(-1)^{\frac x4}\frac vu\pmod p&\text{if }4\mid x\text{ and
}b\ev 1,3\pmod 8,\\ 2(-1)^{\frac x4}\frac vu\pmod p&\text{if }4\mid x\text{
  and }b\ev 5,7\pmod 8.\end{cases}$$

Sun has checked this conjecture for $b<60$ and $p<20000$. 

{\bf Conjecture 2} Let $p\ev 1\pmod 4$ be a prime, $b\in\ZZ$, $b\ev
4\pmod 8$ and $p=x^2+(b^2/4+1)y^2=u^2+v^2\neq b^2/4+1$ with
$x,y,u,v\in\ZZ$. Suppose that $u$ and all odd parts of $x,y,v$ are
$\ev 1\pmod
4$. Then $$U_{\frac{p-1}4}(b,-1)\ev\begin{cases}(-1)^{\frac{b+4}8+\frac v4}\frac
yx\pmod p&\text{if }2\| x,\\ \frac{yv}{xu}\pmod p&\text{if }2\| y,\\
0\pmod p&\text{if }4\mid xy\end{cases}$$
and
$$V_{\frac{p-1}4}(b,-1)\ev\begin{cases}2(-1)^{\frac{y+v}4}\pmod p&\text{if
}4\mid y,\\ 2(-1)^{\frac{b-4}8+\frac x4}\frac vu\pmod p&\text{if }4\mid x,\\
0\pmod p&\text{if }4\nmid xy.\end{cases}$$

This conjeecture was checked for $b<100$ and $p<20000$. 

{\bf Conjecture 3} Let $p\ev 1\pmod 4$ be a prime, $b\in\ZZ$, $8\mid
b$ and $p=x^2+(b^2/4+1)y^2=u^2+v^2\neq b^2/4+1$ with
$x,y,u,v\in\ZZ$. Suppose that $u$ and all odd parts of $x,y,v$ are
$\ev 1\pmod 4$. Then $$U_{\frac{p-1}4}(b,-1)\ev\begin{cases}0\pmod
p&\text{if }4\mid xy,\\ -(-1)^{(\frac b8-1)y}\frac{yv}{xu}\pmod p&\text{if
}4\nmid xy\end{cases}$$ and
$$V_{\frac{p-1}4}(b,-1)\ev\begin{cases}2(-1)^{\frac{xy+v}4+\frac b8y}\pmod
p&\text{if }4\mid xy,\\ 0\pmod p&\text{if }4\nmid xy.\end{cases}$$

Conjecture 3 was checked for $b<100$ and $p<20000$. 

{\bf Conjecture 4} Let $p\ev 1\pmod 4$ be a prime, $b\in\ZZ$, $b\ev
2\pmod 4$ and $p=x^2+(b^2/4+1)y^2=u^2+v^2$ with
$x,y,u,v\in\ZZ$. Suppose that $u$ and the odd parts of $x,y$ are $\ev
1\pmod 4$. Then
$$U_{\frac{p-1}4}(b,-1)\ev\begin{cases}(-1)^{\frac{b+v-2}4}\frac yx\pmod
p&\text{if }2\| y,\\ 0\pmod p&\text{if }4\mid y\end{cases}$$
and
$$V_{\frac{p-1}4}(b,-1)\ev\begin{cases}0\pmod p&\text{if }2\| y,\\
2(-1)^{\frac{y+v}4}\pmod p&\text{if }4\mid y.\end{cases}$$

Conjecture 4 has been checked for $b<100$ and $p<20000$. 

\bff In all four conjectures we have $p=x^2+dy^2=u^2+v^2$, where
$d=b^2+4$ in Conjecture 1, when $b$ is odd, and $d=b^2/4+1$ in
Cojectures 2, 3 and 4, when $b$ is even. We have
$\frac{V_n(b,-1)+U_n(b,-1)\sqrt{b^2+4}}2=\varepsilon^n$, where
$\varepsilon:=\frac{b+\sqrt{b^2+4}}2$. We have $\sqrt{b^2+4}=\sqrt d$ or $2\sqrt d$,
corresponding to $b$ odd or even, respectively. So $\varepsilon =\frac{b+\sqrt d}2$
and $\frac{V_n(b,-1)+U_n(b,-1)\sqrt d}2=\varepsilon^n$ if $b$ is odd and $\varepsilon =\frac
b2+\sqrt d$ and $\frac{V_n(b,-1)}2+U_n(b,-1)\sqrt d=\varepsilon^n$ if $b$ is even.
\eff

\bff We make some remarks regarding the values of $d,p,x,y,u,v$ modulo
powers of $2$. 

Note that if $b\ev 2\pmod 4$, i.e. in the case of Conjecture 4, we
have that $d=b^2/4+1$ is even. More precisely $d\ev 2\pmod 8$. In
Conjecture 3 $d=b^2/4+1\ev 1\pmod 8$. In Conjectures 1 and 2 we have
$d\ev 5\pmod 4$. 

Since $p=u^2+v^2$ we have $p\ev 1\pmod 8$ if $4\mid v$ and $p\ev
5\pmod 8$ if $2\| v$. 

If $d\ev 5\pmod 8$, that is in the cases of Conlectures 1 and 2, we
have $p=x^2+dy^2\ev 1\pmod 8$ if $2\| x$ or $4\mid y$ and $p\ev 5\pmod
8$ if $4\mid x$ or $2\| y$. 

If $d\ev 1\pmod 8$, in the case of Conjecture 3, we have
$p=x^2+dy^2\ev 1\pmod 8$ if $4\mid x$ or $4\mid y$ and $p\ev 5\pmod 8$
if $2\| x$ or $2\| y$. 

If $d\ev 2\pmod 8$, in the case of Conjecture 4, we have
$p=x^2+dy^2\ev 1\pmod 4$, which implies that $y$ is even and so $p\ev
1\pmod 8$, which implies $4\mid v$. 
\eff

\blm If any of the Conjectures 1-4 is true for some $b\in\ZZ$ and some
prime $p$, then it is true for $-b$ and $p$. 
\elm
\pf Let $A,B$ be the predicted values $\mod p$ for
$U_{\frac{p-1}4}(b,-1)$ and $V_{\frac{p-1}4}(b,-1)$, respectively, and let
$A',B'$ be the similar values for $U_{\frac{p-1}4}(-b,-1)$ and
$V_{\frac{p-1}4}(-b,-1)$. Since
$V_{\frac{p-1}4}(-b,-1)=(-1)^{\frac{p-1}4}V_{\frac{p-1}4}(b,-1)$ and
$U_{\frac{p-1}4}(b,-1)=(-1)^{\frac{p+3}4}U_{\frac{p-1}4}(b,-1)$, we have to
prove that $A'=(-1)^{\frac{p+3}4}A$ and $B'=(-1)^{\frac{p-1}4}B$. We
consider the cases that occur. 

For Conjecture 1(i) we have $B=B'=0$ so $B'=(-1)^{\frac{p-1}4}B$ holds
trivially. If $2\| x$ note that one of $b$ and $-b$ is $\ev 1,3\pmod
8$ and the other is $\ev 5,7\pmod 8$. So one of $A,A'$ is $(-1)^{\frac
u4}\frac{2y}x$ and the other one is $-(-1)^{\frac u4}\frac{2y}x$. Thus $A=-A'$
and we have to prove that $(-1)^{\frac{p+3}4}=-1$, i.e. that $\frac{p+3}4$
is odd or, equivalently, $p\ev 1\pmod 8$. By 1.2 this follows from
$2\| x$. If $2\| y$ then $A=A'=\frac{2yu}{xv}$ and we have to prove that
$\frac{p+3}4$ is even, i.e. that $p\ev 5\pmod 8$. This follows from $2\|
y$ by 1.2.

For Conjecture 1(ii) we have $A=A'=0$ so $A'=(-1)^{\frac{p+3}4}A$ holds
trivially. If $4\mid y$ then $B=B'$ so we have to prove that
$(-1)^{\frac{p-1}4}=1$, i.e. that $p\ev 1\pmod 8$. This follows from
1.2. If $4\mid x$ then one of $b$ and $-b$ is $\ev 1,3\pmod 8$ and the
other is $\ev 5,7\pmod 8$ so we get $B=-B'$ and we have to prove that
$(-1)^\frac{p-1}4=-1$, i.e. that $p\ev 5\pmod 8$, which follows from
$4\mid x$ by 1.2.

For Conjecture 2 if $2\| x$ then $A=(-1)^{\frac{b+4}8+\frac u4}\frac vu$ and
$A'=(-1)^{\frac{-b+4}8+\frac u4}\frac vu$ so we have to prove that
$(-1)^{\frac{p+3}4}=A/A'=(-1)^{\frac b4}=-1$, i.e. that $p\ev 1\pmod 8$,
which follows from $2\| x$ by 1.2. If $2\| y$ then $A=A'$ so we have
to prove that $(-1)^{\frac{p+3}4}=1$, i.e. that $p\ev 5\pmod 8$. This
follows from $2\| y$ by 1.2. If $4\mid y$ then $B=B'$ so we have to
prove that $(-1)^{\frac{p-1}4}=1$, i.e. that $p\ev 1\pmod 8$, which
follws from $4\mid y$ by 1.2. If $4\mid x$ then we have to prove that
$(-1)^{\frac{p-1}4}=B/B'=(-1)^{\frac b4}=-1$, i.e. that $p\ev 5\pmod 8$,
which follows from $4\mid x$ by 1.2. Finally, if $4\mid xy$ then
$A'=(-1)^{\frac{p+3}4}A=0$ and if $4\nmid xy$ then
$B'=(-1)^{\frac{p-1}4}B=0$.

For Conjecture 3 if $4\mid xy$ then $A'=(-1)^{\frac{p+3}4}A=0$. If
$4\nmid xy$ then we have to prove that $(-1)^{\frac{p+3}4}=A/A'=(-1)^{\frac
b4y}=1$, i.e. that $p\ev 5\pmod 8$. But this follows from $4\nmid xy$
by 1.2. (We have either $2\| x$ or $2\| y$.) If $4\mid xy$ then we
have to prove that $(-1)^{\frac{p-1}4}=B/B'=(-1)^{\frac b4y}=1$, i.e. that
$p\ev 1\pmod 8$. This follows from 1.2, since we have either $4\mid x$
or $4\mid y$. If $4\nmid xy$ then $B'=(-1)^{\frac{p-1}4}B=0$. 

For Conjecture 4 note that $p\ev 1\pmod 8$, by 1.2. If $2\| y$ then
$A/A'=(-1)^{\frac b2}=-1=(-1)^{\frac{p+3}4}$. If $4\mid y$ then
$A'=(-1)^{\frac{p+3}4}A=0$. If $2\| y$ then $B'=(-1)^{\frac{p-1}4}B=0$. If
$4\mid y$ then $B=B'$ and $(-1)^{\frac{p-1}4}=1$ so we are done. \qed


\section{Proof of the cojectures}

We will use methods from class field theory which are similar to those
from [B].  Given a number field $F$ and a (possibly archimedian)
prime $\p$ of $F$ for any $x\in F$ we denote by $x_\p$ its image in
$F_\p$. When there is no danger of confusion we simply write $x$
instead of $x_\p$. If $E/F$ is a finite abelian extension and $\P$ is
a prime of $E$ staying over $\p$ then we denote by $\legendre{\cdot,E/F}\p
:F_\p^\ti\to Gal(E/F)$ the Artin map and by $(\cdot,E_\P/F_\p
):F_\p^\ti\to Gal(E_\P/F_\p )$ the local Artin map. If we identify $
Gal(E_\P/F_\p )$ with its image in $Gal(E/F)$ then $\legendre{a,E/F}\p
=(a,E_\P/F_\p )$ for any $a\in F_\p^\ti$. 

Denote $F=\QQ (\sqrt d)$ and $E=F(\mu_4)=\QQ (\sqrt d,i)$. Let $A_1:=(x-y\sqrt
di)(u+vi)\1$. Let $L=E(\sqrt[4]{A_1})$. Since $A_1$ is not a square in
$E$ and $\mu_4\sb E$, we have $Gal(L/E)\ap\ZZ_4=\la\s\ra$, where
$\sigma\in Gal(L/E)$ is given by $\sqrt[4]{A_1}\mapsto
i\sqrt[4]{A_1}$.

\blm The extension $L/F$ is Galois and $Gal(L/F)\ap\ZZ_2\ti\ZZ_4$.
\elm
\pf First we prove that the exension $L/F$ is Galois. Let
$\alpha_1=\sqrt[4]{A_1}$. Let $\( F$ be an algebraic closure of $F$ containing
$L$. Let $\alpha\in\( F$ be a conjugate of $\alpha_1$ over $F$. We have to prove that
$\alpha\in L$. Now $\alpha^4$ is a conjugate of $\alpha_1^4=A_1$ over $F$ so $\alpha^4\in\{
A_1,A_3\}$, where $A_3$ is the image of $A_1$ under the automorphism $i\mapsto
i^3=-i$ of $Gal(E/F)$. 

We have $A_1=(x-y\sqrt di)(u+vi)\1$ so $A_3=(x+y\sqrt di)(u-vi)\1$. But $(x+y\sqrt
di)(x-y\sqrt di)=(u+vi)(u-vi)=p$ so $A_3=A_1\1 =\alpha_3^4$, where
$\alpha_3=\alpha_1\1$. Hence $\alpha_3\in L$. 

Now $\alpha^4=A_k=\alpha_k^4$ with $k\in\ZZ_4^\ti$ so $\alpha=i^l\alpha_k$ for some $l$. Hence
the 8 conjugates of $\alpha_1$ over $F$ are $i^l\alpha_k$ with $l\in\ZZ_4$ and
$k\in\ZZ_4^\ti$. They all belong to $L$, since $\alpha_1,\alpha_3,i\in L$. 

Let $\phi\in Gal(L/F)$. Then $\phi_{|E}$ is given by $i\mapsto i^k$ for some
$k\in\ZZ_4^\ti$. It follows that $\phi (\alpha_1)^4=\phi (A_1)=A_k=\alpha_k^4$ so
$\phi (\alpha_1)=i^l\alpha_k$ for some $l$. Thus the elements of $Gal(L/F)$ are given
by $i\mapsto i^k$, $\alpha_1\mapsto i^l\alpha_k$ with $k\in\ZZ_4^\ti$,
$l\in\ZZ_4$. For $k\in\ZZ_4^\ti$ we denote by $\tau_k$ the morphism $i\mapsto
i^k$, $\alpha_1\mapsto\alpha_k$. In particular, $\tau_1=1_L$. 

Let $H=\{\tau_1,\tau_3\}$. We prove that $H$ is a subgroup of $Gal(L/F)$, the
mapping $k\mapsto\tau_k$ is a morphism between $\ZZ_4^\ti$ and $H$, and
$Gal(L/F)$ is the internal product of $H$ and $Gal(L/E)=\la\s\ra$. To do this it
is enough to prove that $\tau_3^2=1_L$ and $\tau_3\s =\s\tau_3$. 
We have $\tau_3(i)=i^3$ and $\tau_3(\alpha_1)=\alpha_3=\alpha_1\1$ so
$\tau_3^2(i)=i^9=i$ and $\tau_3^2(\alpha_1)=\tau_3(\alpha_1\1 )=\alpha_1$. Hence
$\tau_3^2=1_L$. We have $\tau_3\s (i)=\tau_3 (i)=i^3$ and $\s\tau_3(i)=\s
(i^3)=i^3$. Also $\tau_3\s (\alpha_1)=\tau_3(i\alpha_1)=i^3\alpha_3=-i\alpha_1\1$ and
$\s\tau_3(\alpha_1)=\s (\alpha_1\1 )=(i\alpha_1)\1 =-i\alpha_1\1$. So $\tau_3\s =\s\tau_3$. 

In conclusion, $Gal(L/F)\ap H\ti
Gal(L/E)\ap\ZZ_4^\ti\ti\ZZ_4\ap\ZZ_2\ti\ZZ_4$. \qed

We define a character $\chi :Gal(L/F)\z\mu_4$ by $\chi
(\s^k\tau_l)=i^k$. We have
$$\prod_\q\legendre{\varepsilon,L/F}\q =1\text{ so
}\prod_\q\chi\left(\legendre{\varepsilon,L/F}\q\right)=1.$$

Note
that $\s^k(\alpha_1)a_1\1=i^k$ so for any $\phi\in\la\s\ra =Gal(L/E)$ we
have $\chi (\phi )=\phi (\alpha_1)\alpha_1\1$. 

Now $\legendre dp=\legendre{-1}p=1$ so $p$ splits completely in $F=\QQ (\sqrt d)$
and in $E=\QQ (\sqrt d,i)$. Let $\P$ be a prime of $E$ staying over
$p$. We have $p=(x+y\sqrt di)(x-y\sqrt di)=(u+vi)(u-vi)$ and $\od_\P p=1$ so
for exactly one choice of the $\pm$ sign we have $x\pm y\sqrt di\in\P$
and for exactly one choice of the $\pm$ sign we have $u\pm
vi\in\P$. Moreover for these choices of the $\pm$ sign we have $\od_\P
(x\pm y\sqrt di)=\od_\P (u\pm vi)=1$. Of the four primes of $E$ staying
over $p$ we choose $\P$ s.t. $x-y\sqrt di,u-vi\in\P$. Let $\p =\P\cap F$
and let $\p'$ be the other prime of $F$ staying over $p$. 

Denote by $\j_\pm$ the two archimedian primes of $F$ corresponding to
the two embeddings $F\hookrightarrow\RR$ given by $\sqrt d\mapsto\pm\sqrt
d$. 

\blm (i) $\chi\left(\legendre{\varepsilon,L/F}\p\right)\ev\varepsilon^{-\frac{p-1}4}\pmod\P$.

(ii) $\chi\left(\legendre{\varepsilon,L/F}\q\right)=1$ if $\q\neq\p$ and $\q\nmid
2$. 
\elm
\pf We have $\od_\P A_1=\od_\P (x-y\sqrt di)-\od_\P (u+vi)=1-0=1$. 

Let $\phi\in Gal(E/\QQ )$ be given by $\sqrt d\mapsto -\sqrt d$ and
$i\mapsto i$. Note that $\phi^2=1_E$ so $\phi\1 (A_1)=\phi
(A_1)=(x+y\sqrt di)(u+vi)\1$. Let $\P'=\phi (\P )$. Since $\phi_{|F}\neq
1_F$ we have $\phi (\p )=\p'$. We have $\od_{\P'}A_1=\od_\P\phi\1
(A_1)=\od_\P (x+y\sqrt di)-\od_\P (u+vi)=0-0=0$. 

Since $\p$ splits in $E$ we have $[E_\P :F_\p]=1$ so $\N_{E_\P/F_\p}\varepsilon
=\varepsilon$ so $\legendre{\varepsilon,L/F}\p =\legendre{\varepsilon,L/E}\P\in Gal(L/E)$. It follows that
$\chi\left(\legendre{\varepsilon,L/F}\p\right) =\legendre{\varepsilon,L/E}\p (\alpha_1)\alpha_1\1
=\legendre{A_1,\varepsilon}\P_4=\legendre{\varepsilon,A_1}\P_4\1$. (We have $L=E(\alpha_1)$,
$\alpha_1^4=A_1\in E$ and $\mu_4\sb E$.) 

Since $\varepsilon$ is a unit in $F_\p$ we have
$\legendre{\varepsilon,A_1}\P_4\ev\varepsilon^{\frac{\N\P^{\od_\P A_1}-1}4}\pmod\P$. But $\od_\P
A_1=1$ and, since $p$ splits completely in $E$, we have $\N\P
=p$. Hence $\legendre{\varepsilon,A_1}\P_4\ev\varepsilon^{\frac{p-1}4}\pmod\P$ and
$\chi\left(\legendre{\varepsilon,L/F}\p\right)\ev\varepsilon^{-\frac{p-1}4}\pmod\P$. Similarly
$\chi\left(\legendre{\varepsilon,L/F}{\p'}\right) =\legendre{\varepsilon,A_1}{\P'}_4\1$. But both
$\varepsilon,A_1$ are units in $F_{\P'}$ so $\legendre{\varepsilon,A_1}{\P'}_4=1$. 

Hence we have proved (i) and (ii) in the case $\q =\p'$. We now prove
(ii) in the remaining cases. 

Suppose first that $\q$ is non-archimedian, $\q\nmid 2p$. It is enough
to prove that $\legendre{\varepsilon,L/F}\q =1_L$. But $\varepsilon$ is a unit in $F_\q$ so
it is enough to prove that $\q$ doesn't ramify in $L$. Now $\q\nmid 2$
so $\q$ doesn't ramify in $E=F(i)$. Hence it is enough to prove that
$\Q$ doesn't ramify in $L$, where $\Q$ is some prime of $E$ staying
over $\q$. But this follows from the fact that $\Q\nmid 2$,
$L=E(\sqrt[4]{A_1})$ and $A_1=(x-y\sqrt di)(u+vi)\1$ is unit in
$E_\Q$. (Both $x-y\sqrt di,u+vi$ divide $p$, which is a unit in $E_\Q$,
so they are units in $E_\Q$ as well.)

If $\q=\j_+$ then $\varepsilon_\q =\varepsilon >0$ so $\legendre{\varepsilon,L/F}\q =1_L$ and we are
done. 

If $\q=\j_-$ then let $\Q$ be the infinite prime of $E$ staying over
$\q$ corresponding to the embedding of $E=\QQ (\sqrt d,i)$ in $\CC$ given
by $\sqrt d\mapsto -\sqrt d$ and $i\mapsto i$. Also take $\cal Q$ a prime of
$L$ stayiong over $\Q$. 

We have $F_\q\ap\RR$ and $E_\Q\ap L_{\cal Q}\ap\CC$. Now $\varepsilon_\q
=\(\varepsilon$, where $\(\varepsilon$ is the conjugate of $\varepsilon$. But $\varepsilon >0$ and $\varepsilon\(\varepsilon
=\N (\varepsilon )=-1$ so $\(\varepsilon <0$. So if $\phi =\legendre{\varepsilon,L/F}\q$ then $\phi$
corresponds to $c$, the conjugacy automorphism of $Gal(L_{\cal Q}/F_\q
)$. So we want to determine the automorphism $\phi\in Gal(L/F)$
corresponding to $c\in Gal(L_{\cal Q}/F_\q )$. We have $c(i_{\cal
Q})=c(i)=-i=-i_{\cal Q}$ so $\phi (i)=-i=i^3=\tau_3(i)$ so
$\phi_{|E}={\tau_3}_{|E}$. Hence $\phi\in\tau_3
Gal(L/E)=\tau_3\la\s\ra$. So $\phi =\s^k\tau_3$ for some $k$ and we
have $\chi\left(\legendre{\varepsilon,L/F}\q\right) =i^k$. 

Now $\phi (\alpha_1)=\s^k\tau_3 (\alpha_1)=\s^k(\alpha_1\1 )=i^{-k}\alpha_1\1$ so
$c(\alpha_{1{\cal Q}})=i^{-k}\alpha_{1{\cal Q}}\1$ so $\alpha_{1{\cal Q}}c(\alpha_{1{\cal
Q}})=i^{-k}$. But  $\alpha_{1{\cal Q}}c(\alpha_{1{\cal Q}})=\alpha_{1{\cal
Q}}\({\alpha_{1{\cal Q}}}\in\RR_+$. So $i^{-k}\in\RR_+$, which implies
$i^{-k}=1$. Hence $\chi\left(\legendre{\varepsilon,L/F}\q\right) =i^k=1$. \qed

Since $\prod_\q\chi\left(\legendre{\varepsilon,L/F}\q\right) =1$, Lemma 2.2 implies

\bco $\varepsilon^{\frac{p-1}4}\ev\prod_{\q\mid
 2}\chi\left(\legendre{\varepsilon,L/F}\q\right)\pmod\P$. 
\eco

\blm If $\alpha\in\{\pm 1\}$, then: 

(i) We have $\prod_{\q\mid 2}\chi\left(\legendre{\varepsilon,L/F}\q\right) =\alpha$ iff
$p\ev 1\pmod 8$ and $\varepsilon^{\frac{p-1}4}\ev\alpha\pmod p$ or $p\ev 5\pmod 8$ and
$\varepsilon^{\frac{p-1}4}\ev -\alpha\frac{yv}{xu}\sqrt d\pmod p$. 

(ii) We have $\prod_{\q\mid 2}\chi\left(\legendre{\varepsilon,L/F}\q\right) =\alpha i$
iff $p\ev 1\pmod 8$ and $\varepsilon^{\frac{p-1}4}\ev -\alpha\frac yx\sqrt d\pmod p$ or
$p\ev 5\pmod 8$ and $\varepsilon^{\frac{p-1}4}\ev -\alpha\frac vu\pmod p$. 
\elm
\pf Let $\eta =\prod_{\q\mid 2}\chi\left(\legendre{\varepsilon,L/F}\q\right)$. We
have $\eta\in\mu_4$. By Corollary 2.3,
$\varepsilon^{\frac{p-1}4}\ev\eta\pmod\P$. Obviously $\varepsilon^{\frac{p-1}4}\mod p$
uniquely determines $\varepsilon^{\frac{p-1}4}\mod\P =\eta\mod\P$. In fact it
uniquely determines $\eta\in\mu_4$ since if $\eta_1,\eta_2\in\mu_4$
and $\eta_1\ev\eta_2\pmod\P$ then $\eta_1=\eta_2$. (Otherwise
$\eta_1-\eta_2\mid 2$ and so $\eta_1-\eta_2\notin\P$.) 

Conversely we want to show that $\eta$ uniquely determines
$\varepsilon^{\frac{p-1}4}\mod p$, as described in (i) and (ii). We have
$\varepsilon^{\frac{p-1}4}=A+B\sqrt d$ with $2A,2B\in\ZZ$. If $\(\varepsilon$ is the conjugate
of $\varepsilon$ then $\(\varepsilon^{\frac{p-1}4}=A-B\sqrt d$. But $\varepsilon\(\varepsilon =\N\varepsilon =-1$ so
$\(\varepsilon =-\varepsilon\1$. Since $A+B\sqrt d=\varepsilon^{\frac{p-1}4}\ev\eta\pmod\P$ we also
have $A-B\sqrt d=\(\varepsilon^{\frac{p-1}4}=(-1)^{\frac{p-1}4}\varepsilon^{-\frac{p-1}4}\ev
(-1)^{\frac{p-1}4}\eta\1\pmod\P$. Therefore
$2A\ev\eta+(-1)^{\frac{p-1}4}\eta\1\pmod\P$ and 
$2B\sqrt d\ev\eta-(-1)^{\frac{p-1}4}\eta\1\pmod\P$. Hence we know
$2A,2B\mod\P$ in terms of $\eta$. But $2A,2B\in\ZZ$ and $\P\cap\ZZ
=p\ZZ$ so in fact we know $2A,2B\mod p$ and so we know
$\varepsilon^{\frac{p-1}4}\mod p$. We now discuss the cases that occur. 

(i) If $\eta =\alpha$ then $\eta\1 =\alpha$ so
$(-1)^{\frac{p-1}4}\eta\1 =(-1)^{\frac{p-1}4}\alpha$. If $p\ev 1\pmod 8$ then
$(-1)^{\frac{p-1}4}\eta\1 =\alpha$ so we get $2A\ev\alpha +\alpha=2\alpha\pmod\P$ and $2B\sqrt
d\ev\alpha-\alpha =0\pmod\P$. It follows that $A\ev\alpha\pmod p$ and $B\ev 0\pmod p$. So
$\varepsilon^{\frac{p-1}4}=A+B\sqrt d\ev\alpha\pmod p$, as claimed. If $p\ev 5\pmod 8$
then $(-1)^{\frac{p-1}4}\eta\1 =-\alpha$ so we get $2A\ev\alpha +(-\alpha)=0\pmod\P$ and
$2B\sqrt d\ev\alpha-(-\alpha) =2\alpha\pmod\P$. It follows that $A\ev 0\pmod p$. On the other
hand $x\ev y\sqrt di\pmod\P$ and $u\ev vi\pmod\P$ so $2xuB\sqrt d\ev
2xu\alpha\ev -2yv\alpha\sqrt d\pmod\P$, which implies that $2xuB\ev
-2yv\alpha\pmod\P$. But both sides belong to $\ZZ$ so the congruence also
holds $\mod p$. It follows that $B\ev -\alpha\frac{yv}{xu}\pmod p$ and
so $\varepsilon^{\frac{p-1}4}\ev -\alpha\frac{yv}{xu}\sqrt d\pmod p$. 

(ii) If $\eta =\alpha i$ then $\eta\1 =-\alpha i$ so
$(-1)^{\frac{p-1}4}\eta\1 =(-1)^{\frac{p+3}4}\alpha i$. If $p\ev 1\pmod 8$ then
$(-1)^{\frac{p-1}4}\eta\1 =-\alpha i$ so we get $2A\ev\alpha i+(-\alpha i)=0\pmod\P$
and $2B\sqrt d\ev\alpha i-(-\alpha i)=2\alpha i\pmod\P$. It follows that $A\ev 0\pmod
p$. We also have $x\ev y\sqrt di\pmod\P$ and so $2xB\sqrt d\ev 2x\alpha i\ev
-2y\alpha\sqrt d\pmod\P$, so $2xB\ev -2y\alpha\pmod\P$. This implies
that $2xB\ev -2y\alpha\pmod p$ so $B\ev -\alpha\frac yx\pmod p$. Hence
$\varepsilon^{\frac{p-1}4}\ev -\alpha\frac
yx\sqrt d\pmod p$. If $p\ev 5\pmod 8$ then $(-1)^{\frac{p-1}4}\eta\1 =\alpha i$
so we get $2A\ev\alpha i+\alpha i=2\alpha i\pmod\P$ and $2B\sqrt d\ev\alpha i-\alpha
i=0\pmod\P$. It follows that $B\ev 0\pmod p$. Since also $u\ev
vi\pmod\P$, we get $2uA\ev 2u\alpha i\ev -2v\alpha\pmod\P$. It follows that
$2uA\ev -2v\alpha\pmod p$ and so $A\ev -\alpha\frac vu\pmod p$. Hence
$\varepsilon^{\frac{p-1}4}\ev -\alpha\frac vu\pmod p$. \qed  

\bco Let $\eta =\prod_{\q\mid 2}\chi\left(\legendre{\varepsilon,K/F}\q\right)$. Then

(i) Conjecture 1 is equivalent to
$$\eta =\begin{cases}-(-1)^{\frac v4}i&\text{if }2\| x\text{ and }b\ev
1,3\pmod 8\\
(-1)^{\frac v4}i&\text{if }2\| x\text{ and }b\ev 5,7\pmod 8\\
-1&\text{if }2\| y\\
(-1)^{\frac{y+v}4}&\text{if }4\mid y\\
(-1)^{\frac x4}i&\text{if }4\mid x\text{ and }b\ev 1,3\pmod 8\\
-(-1)^{\frac x4}i&\text{if }4\mid x\text{ and }b\ev 5,7\pmod 8.
\end{cases}$$

(ii) Conjecture 2 is equivalent to
$$\eta =\begin{cases}-(-1)^{\frac{b+4}8+\frac v4}i&\text{if }2\| x\\
-1&\text{if }2\| y\\
(-1)^{\frac{y+v}4}&\text{if }4\mid y\\
(-1)^{\frac{b+4}8+\frac x4}i&\text{if }4\mid x.
\end{cases}$$

(iii) Conjecture 3 is equivalent to
$$\eta =\begin{cases}(-1)^{(\frac b8-1)y}&\text{if }4\nmid xy\\
(-1)^{\frac{xy+v}4+\frac b8y}&\text{if }4\mid xy.
\end{cases}$$

(iv) Conjecture 4 is equivalent to
$$\eta =\begin{cases}(-1)^{\frac{b+v+2}4}i&\text{if }2\|y\\
(-1)^{\frac{y+v}4}&\text{if }4\mid y.
\end{cases}$$
\eco
\pf As seen in 1.1 in the case of Conjecture 1 we have
$\frac{V_{\frac{p-1}4}(b,-1)+U_{\frac{p-1}4}(b,-1)\sqrt
d}2=\varepsilon^{\frac{\p-1}4}$. Then, by Lemma 2.4, for $\alpha\in\{\pm 1\}$
we have: 

If $p\ev 1\pmod 8$ then $\eta =\alpha$ is equivalent to
$V_{\frac{p-1}4}(b,-1)\ev 2\alpha \pmod p$ and $U_{\frac{p-1}4}(b,-1)\ev 0\pmod
p$; and $\eta =\alpha i$ is equivalent to $V_{\frac{p-1}4}(b,-1)\ev 0$ and
$U_{\frac{p-1}4}(b,-1)\ev -2\alpha\frac yx\pmod p$. 

If $p\ev 5\pmod 8$ then $\eta =\alpha$ is equivalent to $V_{\frac{p-1}4}(b,-1)\ev
0\pmod p$ and $U_{\frac{p-1}4}(b,-1)\ev -2\alpha\frac{yv}{xu}\pmod p$;
and $\eta =\alpha i$
is equivalent to $V_{\frac{p-1}4}(b,-1)\ev -2\alpha\frac vu\pmod p$ and
$U_{\frac{p-1}4}(b,-1)\ev 0\pmod p$. 

In the cases of Conjectures 2, 3 and 4 we have
$\frac{V_{\frac{p-1}4}(b,-1)}2+U_{\frac{p-1}4}(b,-1)\sqrt d=\varepsilon^{\frac{\p-1}4}$. So
the factor $2$ in the formulas above for $U_{\frac{p-1}4}(b,-1)\mod p$
should be dropped. That is, $-2\alpha\frac yx$ and $-2\alpha\frac{yv}{xu}$ should be
replaced by $-\alpha\frac yx$ and $-\alpha\frac{yv}{xu}$, respectively. 

We prove (i). If $2\| x$ then Conjecture 1 states that
$V_{\frac{p-1}4}(b,-1)\ev 0\pmod p$ and $U_{\frac{p-1}4}(b,-1)\ev (-1)^{\frac
  v4}2\frac yx\pmod p$ or $-(-1)^{\frac v4}2\frac  yx\pmod p$, corresponding to
$b\ev 1,3\pmod 8$ or $b\ev 5,7\pmod 8$, respectively. But from 1.2 we
have $p\ev 1\pmod 8$ so this is equivalent to $\eta =-(-1)^{\frac v4}i$
or $\eta =(-1)^{\frac v4}i$, respectively. If $2\| y$ then Conjecture 1
states that $V_{\frac{p-1}4}(b,-1)\ev 0\pmod p$ and
$U_{\frac{p-1}4}(b,-1)\ev 2\frac{yv}{xu}\pmod p$. Since by 1.2 $p\ev 5\pmod
8$ this is equivalent to $\eta =-1$. If $4\mid y$ then Conjecture 1
states that $V_{\frac{p-1}4}(b,-1)\ev 2(-1)^{\frac{y+v}4}\pmod p$ and
$U_{\frac{p-1}4}(b,-1)\ev 0\pmod p$, which is equivalent to $\eta
=(-1)^{\frac{y+v}4}$. (We have $p\ev 1\pmod 8$.) If $4\mid x$ then
Conjecture 1 states that $V_{\frac{p-1}4}(b,-1)\ev -2(-1)^{\frac x4}\frac vu$
or $2(-1)^{\frac x4}\frac vu\pmod p$, corresponding to $b\ev 1,3\pmod 8$ or
$b\ev 5,7\pmod 8$, respectively, and $U_{\frac{p-1}4}(b,-1)\ev 0\pmod
p$. But this is equivalent to $\eta =(-1)^{\frac x4}$ or $-(-1)^{\frac x4}$,
respectively. (We have $p\ev 5\pmod 8$.) 

We prove (ii). If $2\| x$ then Conjecture 2 states that
$V_{\frac{p-1}4}(b,-1)\ev 0\pmod p$ and $U_{\frac{p-1}4}(b,-1)\ev
(-1)^{\frac{b+4}8+\frac v4}\frac yx\pmod p$. Since $p\ev 1\pmod 8$ this is
equivalent to $\eta =-(-1)^{\frac{b+4}8+\frac v4}i$. If $2\| y$ then
Conjecture states that $V_{\frac{p-1}4}(b,-1)\ev 0\pmod p$ and
$U_{\frac{p-1}4}(b,-1)\ev\frac{yv}{xu}\pmod p$. Since $p\ev 5\pmod 8$ this
is equivalent to $\eta =-1$. If $4\mid y$ then Conjecture 2 states
that $V_{\frac{p-1}4}(b,-1)\ev 2(-1)^{\frac{y+v}4}\pmod p$ and
$U_{\frac{p-1}4}(b,-1)\ev 0\pmod p$. Since $p\ev 1\pmod 8$ this is
equivalent to $\eta =(-1)^{\frac{y+v}4}$. If $4\mid x$ then Conjecture 2
states that $V_{\frac{p-1}4}(b,-1)\ev 2(-1)^{\frac{b-4}8+\frac x4}\frac vu\pmod p$
and $U_{\frac{p-1}4}(b,-1)\ev 0\pmod p$ which is equivalent to $\eta
=-(-1)^{\frac{b-4}8+\frac x4}i=(-1)^{\frac{b+4}8+\frac x4}i$ since $p\ev 5\pmod 8$. 

We prove (iii). By 1.2 if $4\mid xy$ then $p\ev 1\pmod 8$, while if
$2\| xy$ then $p\ev 5\pmod 8$. If $4\mid xy$ then Conjecture 3 states
that $V_{\frac{p-1}4}(b,-1)\ev 2(-1)^{\frac{xy+v}4+\frac b8y}\pmod p$ and
$U_{\frac{p-1}4}(b,-1)\ev 0\pmod p$. Since $p\ev 1\pmod 8$ this is
equivalent to $\eta =(-1)^{\frac{xy+v}4+\frac b8y}$. If $2\| xy$ then
Conjecture 3 states that $V_{\frac{p-1}4}(b,-1)\ev 0\pmod p$ and
$U_{\frac{p-1}4}(b,-1)\ev -(-1)^{(\frac b8-1)y}\frac{yv}{xu}\pmod p$. Since
$p\ev 5\pmod 8$ this is equivalent to $\eta =(-1)^{(\frac b8-1)y}$.

We prove (iv). By 1.2 we have $p\ev 1\pmod 8$. If $2\| y$ then
Conjecture 4 states that $V_{\frac{p-1}4}(b,-1)\ev 0\pmod p$ and
$U_{\frac{p-1}4}(b,-1)\ev (-1)^{\frac{b+v-2}4}\frac yx\pmod p$. But $p\ev
1\pmod 8$ so this is equivalent to $\eta
=-(-1)^{\frac{b+v-2}4}i=(-1)^{\frac{b+v+2}4}i$. If $4\mid y$ then Conjecture
4 states that $V_{\frac{p-1}4}(b,-1)\ev 2(-1)^{\frac{y+v}4}\pmod p$ and
$U_{\frac{p-1}4}(b,-1)\ev 0\pmod p$. Since $p\ev 1\pmod 8$ this is
equivalent to $\eta =(-1)^{\frac{y+v}4}$. \qed

Unfortunately the primes of $F$ staying over $2$ are ramified in $L$
so it is very difficult to calculate the Artin symbols
$\legendre{\varepsilon,L/F}\q$ for $\q\mid 2$. To overcame this problem we will show
that $\eta =\prod_{\q\mid 2}\chi\left(\legendre{\varepsilon,L/F}\q\right)$ only
depends on the congruence classes of $b,p,x,y,u,v\mod$ modulo
appropriate powers of $2$. Obviously same happens with the values of
$\eta$ predicted by Corollary 2.5. This way we reduce our conjectures
to verifying them for a finite number of pairs $b,p$ with $b\in\ZZ$
and $p$ prime. 

Let now $b'\in\ZZ$ and $p'$ is a prime that satisfy this satify the
hypothesis of one of the conjectures and let
$d',x',y',u',v',F',E',L',\s',\tau_k',\chi',\eta'$ be the
$d,x,y,u,v,\ldots$ corresponding to $b',p'$. 

\bff We will assume that one of the following happens: 

(1) $b,b'$ are odd and $b\ev b'\pmod 8$. 

(2) $2\| b,b'$ and $b\ev b'\pmod{32}$.

(3) $2^2\| b,b'$ and $b\ev b'\pmod{16}$.

(4) $8\mid b,b'$ and $b\ev b'\pmod{16}$.

Note that the four cases above correspond to Conjectures 1, 4, 2 and
3, respeectively. 

In the cases (1) and (3) we have $d\ev d'\ev 5\pmod 8$ and so $2$ is
inert in $F,F'$. In the case (2) we have $2\| d$ so $2$ ramifies in
$F$. Same for $F'$. Therefore we will denote by $\q,\q'$ the only
primes of $F,F'$ staying over $2$. 

In the case (1) $b,b'$ are odd and $\ev\mod 8$, which implies that
$b^2\ev b'^2\pmod{16}$. Thus $d=b^2+4$ and $d'=b'^2+4$ are odd and
$\ev\mod 16$ and so $d'/d\ev 1\pmod{16}$. In the case (2) $b/2,b'/2$
are odd and $\ev\mod 16$ so $b^2/4\ev b'^2/4\pmod{32}$, which implies
that $d\ev d'\pmod{32}$. Since also $2\| d$, we have $d'/d\ev
1\pmod{16}$. In the case (3) $b/4,b'/4$ are odd and so $b^2/16\ev
1\pmod 8$ which implies that $d=b^2/4+1\ev 5\pmod{32}$. Similarly
$d'\ev 5\pmod{32}$ and we get $d'/d\ev 1\pmod{32}$.

By [B, Lemma 1.6(i)] in all three cases we have $d'/d=\nu^2$ for some
$\nu\in\QQ_2$. Moreover, in the cases (1) and (2) we have $\nu\ev
1\pmod 8$, while in the case (3) $\nu\ev 1\pmod{16}$. Since $d'=\nu^2
d$ we have $\QQ (\sqrt d)_\q\ap\QQ (\sqrt{d'})_{\q'}$ and the mapping $\psi
:F'_{\q'}\z F_\q$ given by $\sqrt{d'}\mapsto\nu\sqrt d$ is an isomorphism of
local fields.

In the case (4) we have $d\ev d'\ev 1\pmod{16}$ so $2$ splits in
$F,F'$. By [B, Lemma 1.6(i)] we have $d=\nu^2$ for some $\nu\in\QQ_2$,
$\nu\ev 1\pmod 8$. We have two primes of $F$ staying over $2$, $\q_1$
and $\q_2$, and they correspond to the embeddings of $F$ in $\QQ_2$
given by $\sqrt d\mapsto\nu$ and $\sqrt d\mapsto -\nu$. We define similarly
$\nu',\q_1',\q_2'$ for $F'$. For $j=1,2$ we have $F_{\q_1}\ap
F'_{\q'_1}\ap\QQ_2$. Let $\psi_j: F'_{\q'_j}\z F_{\q_j}$ be the only
isomorphism between the two local fields. 
\eff
\blm With the settings of 2.6 we have:

(i) If $8\nmid b,b'$ then $\psi (\varepsilon')/\varepsilon =A^4$ for some $A\in F_\q$. 

(ii) If $8\mid b,b'$ then

$$\prod_{\q\mid 2}\chi\left(\legendre{\varepsilon,L/F}\q\right)
=\begin{cases}\chi\left(\legendre{-1,L/F}{\q_2}\right)&\text{ if }16\mid
b\\
\chi\left(\legendre{5,L/F}{\q_1}\right)
\chi\left(\legendre{3,L/F}{q_2}\right)&\text{ if }16\nmid b.\end{cases}$$

Same for $b'$.
\elm
\pf (i) We consider the cases (1)-(3) of 2.6. Recall that in the cases
(1) and (3) $d\ev 5\pmod 8$ and $2$ is inert in $F$, while in the case
(2) $d\ev 2\pmod 8$ and $2$ ramifies in $F$. 

In the case (1) $b,b'$ are odd and $\ev\mod 8$. We have two cases:

a. If $b\ev b'\pmod{16}$ then recall that $\varepsilon'$ is a root of
$X^2-b'X-1$ and so is $\psi (\varepsilon')$. Also $\varepsilon,\(\varepsilon$ are the roots of
$X^2-bX-1$ so $(X-\varepsilon )(X-\(\varepsilon )=X^2-bX-1$. It follows that $(\psi (\varepsilon'
)-\varepsilon )(\psi (\varepsilon')-\(\varepsilon )=\psi (\varepsilon')^2-b\psi (\varepsilon')-1=(b'-b)\psi (\varepsilon')$,
which implies that $\od_2(\psi (\varepsilon' )-\varepsilon )+\od_2(\psi (\varepsilon')-\(\varepsilon
)=\od_2(b'-b)+\od_2\psi (\varepsilon')=\od_2(b'-b)\geq 4$. (Note that $\varepsilon'$ is
a unit in $F'$, so in $F'_{\q'}$. Thus $\psi (\varepsilon')$ is a unit in
$F_\q$.) Now $\varepsilon'=\frac{b'+\sqrt{d'}}2$ so $\psi (\varepsilon')=\frac{b'+\nu\sqrt d}2$ and
$\(\varepsilon =\frac{b-\sqrt d}2$. Thus $\psi (\varepsilon')-\(\varepsilon =\frac{b'-b}2+\frac{\nu +1}2\sqrt
d$. Since $\frac{b'-b}2$ is even and $\nu\ev 1\pmod 8$, so $\frac{\nu +1}2$
is odd, we have $\od_2(\psi (\varepsilon')-\(\varepsilon )=0$ and so $\od_2(\psi
(\varepsilon')-\varepsilon )\geq 4$. But $\varepsilon $ is a unit in $F_\q$ and so $\od (\psi
(\varepsilon')/\varepsilon -1)\geq 4$, i.e. $\psi (\varepsilon')/\varepsilon\ev 1\pmod{16}$. By Lemma
2.7(i) we get $\psi (\varepsilon')/\varepsilon =A^4$ for some $A\in F_\q$.

b. If $b'\ev b+8\pmod{16}$ then we note that $d^2\varepsilon,d^2\(\varepsilon$ are the
roots of $X^2-bd^2X-d^4$. Therefore $(\psi (\varepsilon' )-d^2\varepsilon )(\psi
(\varepsilon')-d^2\(\varepsilon )=\psi (\varepsilon')^2-bd^2\psi (\varepsilon')-d^4=(b'-bd^2)\psi
(\varepsilon')+1-d^4$. We have $d\ev 5\pmod 8$, i.e. $2^2\|d-1$, which implies
that $2^3\| d^2-1$ and $2^4\| d^4-1$. So $d^2\ev 9\pmod{16}$, which
implies $bd^2\ev 9b\ev b+8\ev b'\pmod{16}$. ($b$ is odd so $8b\ev
8\pmod{16}$.) Since also $d^4\ev 1\pmod{16}$, we have $16\mid
(b'-bd^2)\psi (\varepsilon')+1-d^4$ and so $\od_2(\psi (\varepsilon'
)-d^2\varepsilon )+\od_2(\psi (\varepsilon')-d^2\(\varepsilon )\geq 4$. Now $\psi (\varepsilon')-d^2\(\varepsilon
=\frac{b'-d^2b}2+\frac{\nu +d^2}2\sqrt d$. But $bd^2\ev b'\cdot 1=b'\pmod 4$ so
$\frac{b'-bd^2}2$ is even and $\nu\ev d^2\ev 1\pmod 4$ so $\frac{\nu +d^2}2$
is odd. It follows that $\od_2(\psi (\varepsilon')-d^2\(\varepsilon )=0$ so $\od_2(\psi
(\varepsilon')-d^2\varepsilon )\geq 4$. By the same reasoning as in case a. we have
$\psi (\varepsilon')/(d^2\varepsilon )=A'^4$ for some $A'\in F_\q$ so $\psi
(\varepsilon')/\varepsilon=A^4$, where $A=A'\sqrt d\in F_\q$. 

In the case (3) we proceed similarly as for (1). This time $2^2\|
b,b'$ and $b\ev b'\pmod{16}$ and we consider two cases:

a. If $b\ev b'\pmod{32}$ then we use again the relation $\od_2(\psi
(\varepsilon' )-\varepsilon )+\od_2(\psi (\varepsilon')-\(\varepsilon
)=\od_2(b'-b)\geq 5$. This time however $\varepsilon'=\frac{b'}2+\sqrt
d'$ so $\psi (\varepsilon')=\frac{b'}2+\nu\sqrt d$ and $\(\varepsilon
=\frac b2-\sqrt d$. So $\psi (\varepsilon')-\(\varepsilon
=\frac{b'-b}2+(\nu +1)\sqrt d$. But $4\mid\frac{b'-b}2$ and $2\|\nu
+1$ so $\od_2(\psi (\varepsilon')-\(\varepsilon )=1$ and we
get $\od_2(\psi (\varepsilon')-\varepsilon )\geq 4$. From here the proof follows as for
case a. of (1). 

b. If $b'\ev b+16\pmod{32}$ then we use the relation $(\psi (\varepsilon'
)-d^2\varepsilon )(\psi (\varepsilon')-d^2\(\varepsilon )=(b'-bd^2)\psi (\varepsilon')+1-d^4$. Since
$4\mid b$ and $d^2\ev 1\pmod 8$ we have $bd^2\ev b\pmod{32}$ so
$b'-bd^2\ev b'-b\ev 16\pmod{32}$. But, as seen in the case b. of (1),
we also have $2^4\| d^4-1$ so $1-d^4\ev 16\pmod{32}$. In conclusion,
$(b'-bd^2)\psi (\varepsilon')+1-d^4\ev 16(\psi (\varepsilon')+1)\pmod{32}$. But $\psi
(\varepsilon')+1=\frac{b'}2+1+\nu\sqrt d$. Since $\frac{b'}2+1\ev\nu\ev 1\pmod 2$ we
have $\psi (\varepsilon')+1\ev 1+\sqrt d\ev0\pmod 2$ and so $(b'-bd^2)\psi
(\varepsilon')+1-d^4\ev 0\pmod{32}$. It follows that $\od_2(\psi (\varepsilon' )-d^2\varepsilon
)+\od_2(\psi (\varepsilon')-d^2\(\varepsilon )=\od_2((b'-bd^2)\psi (\varepsilon')+1-d^4)\geq
5$. We have $\psi (\varepsilon')-d^2\(\varepsilon =\frac{b'-bd^2}2+(\nu +d^2)\sqrt d$. But
$b'-bd^2\ev 16\pmod{32}$ so $4\mid\frac{b'-bd^2}2$ and $\nu\ev d^2\ev
1\pmod 4$ so $2\|\nu +d^2$ so $\od_2(\psi (\varepsilon')-d^2\(\varepsilon )=1$. It
follows that $\od_2(\psi (\varepsilon')-d^2\varepsilon )\geq 4$ and the proof follows
same as in the case b. of (1).  

In the case of (2) let $\oo_\q$ be the ring of integers in $F_\q$ and
by $\tilde\q$ the completion of $\q$ in $F_\q$, $\tilde\q =\q\oo_\q
=\sqrt d\oo_\q$. Since $2$ ramifies in $F$ we have
$\od_{\tilde\q}2=2$. We use again the relation $\od_{\tilde\q}(\psi
(\varepsilon' )-\varepsilon )+\od_{\tilde\q}(\psi (\varepsilon')-\(\varepsilon )=\od_{\tilde\q}(b'-b)\geq
10$. (We have $32\mid b'-b$ and $\od_{\tilde\q}32=10$.) But $\psi
(\varepsilon')-\(\varepsilon =\frac{b'-b}2+(\nu +1)\sqrt d$. Since $4\mid\frac{b'-b}2$ and
$2\|\nu +1$ we have $\od_{\tilde\q}(\psi (\varepsilon')-\(\varepsilon )=3$. Hence
$\od_{\tilde\q}(\psi (\varepsilon')-\varepsilon )\geq 7$ i.e. $\psi (\varepsilon')-\varepsilon\in
8\tilde\q$. Since also $\varepsilon\in\oo_\q^\ti$ we get $\psi (\varepsilon')/\varepsilon -1\in
8\tilde\q$ so, by [B, Lemma 1.6(i)], $\psi (\varepsilon')/\varepsilon =A^4$ for some
$A\in F_\q$.  

(ii) We have that $b,b'$ are in the case (4) of 2.6 so $2$ splits in
$F$. For $j=1,2$ we have $F_{\q_j}\ap\QQ_2$ and if $\varepsilon_{q_j}$ is the
image of $\varepsilon$ in $F_{\q_j}\ap\QQ_2$ then
$\varepsilon_{\q_1}\varepsilon_{q_2}=\N_{F/Q}(\varepsilon )=-1$. 

Also note that if $j\in\{ 1,2\}$ then $\varepsilon_{\q_j}$ is a unit in
$F_{\q_j}\ap\QQ_2$ so there is some odd $s\in\ZZ$ with $\varepsilon_{\q_j}\ev
s\pmod{16}$. It follows that $\varepsilon_{\q_j}/s\ev 1\pmod{16}$ so, by [B,
Lemma 1.6(i)], $\varepsilon_{\q_j}/s=A^4$ for some $A\in\QQ_2$. We claim that
$\legendre{\varepsilon,L/F}{\q_j}=\legendre{s,L/F}{\q_j}$. This is equivalent to
$\legendre{\varepsilon/s,L/F}{\q_j}=1$. But $\legendre{\varepsilon/s,L/F}{\q_j}$ is the image of
$(\varepsilon/s,L_\Q/F_{\q_j})=(A^4,L_\Q/F_{\q_j})\in Gal(L_\Q/F_{\q_j})$
in $Gal(L/F)$, where $\Q$ is a prime of $L$ staying over $\q_j$. But
$(A^4,L_\Q/F_{\q_j})=(A,L_\Q/F_{\q_j})^4=1$, since
$Gal(L/F)\ap\ZZ_2\ti\ZZ_4$. 

We have $\varepsilon =\frac b2+\sqrt d$ so $\varepsilon_{\q_1}=\frac b2+\nu$. But $\nu\ev
1\pmod 8$ so if $16\mid b$ then $\varepsilon_{\q_1}\ev 1\pmod 8$, while if
$16\nmid b$ then $\varepsilon_{\q_1}\ev 5\pmod 8$. We consider the two cases. 

If $16\mid b$ then $\varepsilon_{\q_1}\ev 1\pmod 8$ and $\varepsilon_{\q_1}\varepsilon_{\q_2}=-1$
so we have either $\varepsilon_{\q_1}\ev 1\pmod{16}$ and $\varepsilon_{\q_2}\ev
-1\pmod{16}$ or $\varepsilon_{\q_1}\ev 9\ev 5^2\pmod{16}$ and $\varepsilon_{\q_2}\ev
7\ev -5^2\pmod{16}$. In the 
first case we have $\legendre{\varepsilon,L/F}{\q_1}=\legendre{1,L/F}{\q_1}=1$ and
$\legendre{\varepsilon,L/F}{\q_2}=\legendre{-1,L/F}{\q_2}$ so $\prod_{\q\mid
2}\chi\left(\legendre{\varepsilon,L/F}\q\right)
=\chi\left(\legendre{-1,L/F}{q_2}\right)$, as claimed. In the second case
we have $\prod_{\q\mid 2}\chi\left(\legendre{\varepsilon,L/F}\q\right)
=\chi\left(\legendre{5^2,L/F}{\q_1}\right)
\chi\left(\legendre{-5^2,L/F}{\q_2}\right)
=\chi\left(\legendre{-1,L/F}{\q_2}\right)\prod_{\q\mid
2}\chi\left(\legendre{5,L/F}\q\right)^2$. So in this 
case we have to prove that $\prod_{\q\mid
2}\chi\left(\legendre{5,L/F}\q\right)^2=1$.

If $16\nmid b$ then $\varepsilon_{\q_1}\ev 5\pmod 8$ and
$\varepsilon_{\q_1}\varepsilon_{\q_2}=-1$ so we have either $\varepsilon_{\q_1}\ev 5\pmod{16}$
and $\varepsilon_{\q_1}\ev 3\pmod{16}$ or $\varepsilon_{\q_1}\ev 13\ev 5\cdot
5^2\pmod{16}$ and $\varepsilon_{\q_1}\ev 11\ev 3\cdot 5^2\pmod{16} $. In the
first case we have $\prod_{\q\mid 2}\chi\left(\legendre{\varepsilon,L/F}\q\right)
=\chi\left(\legendre{5,L/F}{\q_1}\right)\chi\left(\legendre{3,L/F}{\q_2}\right)$
and we are done. In the second case $\prod_{\q\mid
2}\chi\left(\legendre{\varepsilon,L/F}\q\right) =\chi\left(\legendre{5\cdot
5^2,L/F}{\q_1}\right)\chi\left(\legendre{3\cdot 5^2,L/F}{\q_2}\right)
=\chi\left(\legendre{5,L/F}{\q_1}\right)
\chi\left(\legendre{3,L/F}{\q_2}\right)\prod_{\q\mid
2}\chi\left(\legendre{5,L/F}\q\right)^2$ so again we have to prove that
$\prod_{\q\mid 2}\chi\left(\legendre{5,L/F}\q\right)^2=1$. 

So we have reduced (ii) to proving that $\prod_{\q\mid
2}\chi\left(\legendre{5,L/F}\q\right)^2=1$. Note that the primes $\q_j$ of
$F$ ramify in $E=F(i)$. Let $\Q_j$ be the prime of $E$ staying over
$\q_j$. We have $5=\N_{E/F}(1+2i)$ so if $\phi_j:=\legendre{5,L/F}{\q_j}$
then $\phi_j=\legendre{1+2i,L/E}{\Q_j}$. Since $\phi_j\in
Gal(L/E)$ we have $\chi (\phi_j)=\phi_j(\alpha_1)/\alpha_1$. But
$\alpha_1^4=A_1\in E$, $\mu_4\sb E$ and $\phi_j=\legendre{1+2i,L/E}{\Q_j}$ so
$\phi_j(\alpha_1)/a_1=\legendre{A_1,1+2i}{\Q_j}_4$. Hence
$\chi\left(\legendre{5,L/F}{q_j}\right) =\legendre{1+2i,A_1}{\Q_j}_4$, which
implies that
$\chi\left(\legendre{5,K/F}{\q_j}\right)^2=\legendre{1+2i,A_1}{\Q_j}$. So
we have to prove that
$\legendre{1+2i,A_1}{\Q_1}\legendre{1+2i,A_1}{\Q_2}=1$.

Since $F_{\q_j}\ap\QQ_2$ we have $E_{\Q_j}\ap\QQ_2(\sqrt{-1})\ap\QQ
(i)_{(1+i)}$. (Here $(1+i)$ is the ideal generated by $1+i$ and is the
only prime ideal of $\QQ (i)$ staying over $2$.) Recall that the
isomorphisms $F_{\q_j}\z\QQ_2$ are obtained by extending by continuity
the embeddings $F\z\QQ_2$ given by $\sqrt d\mapsto\pm\nu$. The signs $+$
and $-$ correspond to $j=1$ and $2$, respectively. These isomorphisms
can be extended to isomorphisms $E_{\Q_j}\z\QQ (i)_{(1+i)}$ if we
extend by continuity the embeddings $E\z\QQ (i)_{(1+i)}$ given by $\sqrt
d\mapsto\pm\nu$ and $i\mapsto i$. For $\alpha\in E$ we denote by
$\alpha_{\Q_j}$ its image through the isomorphism $F_{\Q_j}\z\QQ
(i)_{(1+i)}$ desribed above. We have $(1+2i)_{\Q_j}=1+2i$ for $j=1,2$
and, since $A_1=(x-y\sqrt d i)(u+vi)\1$ we have ${A_1}_{\Q_1}=(x-y\nu
i)(u+vi)\1$ and ${A_1}_{\Q_2}=(x+y\nu i)(u+vi)\1$. It follows that
$\legendre{1+2i,A_1}{\Q_2}\legendre{1+2i,A_1}{\Q_2}=\legendre{1+2i,(x-y\nu
i)(u+vi)\1}{(1+i)}\legendre{1+2i,(x+y\nu i)(u+vi)\1}{(1+i)}
=\legendre{1+2i,p(u+vi)^{-2}}{(1+i)}=\legendre{p,1+2i}{(1+i)}$. (Note that
$(x-y\nu i)(x+y\nu i)=x^2+y^2\nu^2=x^2+dy^2=p$.) But $p\in\QQ_2$ 
and $\N_{\QQ (i)_{(1+i)}/\QQ_2}(1+2i)=5$ so
$\legendre{p,1+2i}{(1+i)}=\legendre{p,5}2=1$. \qed

\bff We now make the assumption that $p',x',y',u',v'$ are close in the
$2$-adic topology to $p,x,y,u,v$, as follows: 

%

(1) If $d\ev d'\ev 5\pmod 8$, i.e. in the case of Conjectures 1 and 2, we will
assume that one of the following happens:

a. $2\| x$, $2\| x'$ and $\frac v4\ev\frac{v'}4\pmod 2$. 

b. $2\| y$, $2\| y'$. 

c. $4\mid y$, $4\mid y'$ and $\frac{y+v}4\ev\frac{y'+v'}4\pmod 2$. 

d. $4\mid x$, $4\mid x'$ and $\frac x4\ev\frac{x'}4\pmod 2$. 

(2) If $d\ev d'\ev 2\pmod 8$, i.e. in the case of Conjecture 4, we will
assume that one of the following happens:

a. $2\| y$, $2\| y'$ and $\frac v4\ev\frac{v'}4\pmod 2$. 

b. $4\mid y$, $4\mid y'$ and $\frac v4\ev\frac{v'}4\pmod 2$. 

(3) If $d\ev d'\ev 1\pmod 8$, i.e. in the case of Conjecture 3, we will assume
that one of the following happens: 

a. $2\| x$, $2\| x'$.

b. $2\| y$, $2\| y'$. 

c. $4\mid x$, $4\mid x'$ and $\frac{x+v}4\ev\frac{x'+v'}4\pmod 4$. 

d. $4\mid y$, $4\mid y'$ and $\frac{y+v}4\ev\frac{y'+v'}4\pmod 4$.  

Note that in the cases $d\ev d'\ev 1,5\pmod 8$ one of the $2\| x,x'$,
$2\| y,y'$, $4\mid x,x'$ and $4\mid y,y'$ holds so by 1.2 we have
$p\ev p'\pmod 8$, which in turn implies $v\ev v'\pmod 4$. Same happens
if $d\ev d'\ev 2\pmod 8$, when $p\ev p'\ev 1\pmod 8$ and $4\mid
v,v'$.
\eff

\bff Recall that if $d\ev d'\ev 2,5\pmod 8$, since $d$ is not a square
in $\QQ_2$, $2$ doesn't split in $F=\QQ (\sqrt d)$. But also $-1,-d$ are
nonsquares in $\QQ_2$ so $2$ doesn't split in $E=\QQ (\sqrt d,i)$. Let
$\Q$ be the only prime of $E$ staying over $2$ and also over
$\q$. Similarly we define the prime $\Q'$ of $E'$. In 2.6 we have
defined an isomorphism $\psi :F'_{\q'}\z F_{\q }$ by
$\sqrt{d'}\mapsto\nu\sqrt d$, where $\nu\in\QQ_2$, $\nu\ev 1\pmod 8$. Since
$E=F(i)$, $E'=F'(i)$ we can extend $\psi$ to $\psi':E'_{\Q'}\z E_\Q$
by $\sqrt{d'}\mapsto\nu\sqrt d$ and $i\mapsto i$. We extend further $\Q,\Q'$
to primes $\cal Q,\cal Q'$ of $L,L'$ staying over $\Q,\Q'$. 

In the case of Conjecture 3 we have the two primes $\q_1,\q_2$ of
$F$. They extend to the primes $\Q_1,\Q_2$ of $E$. We have the
isomorphisms $F_{\q_j}\z\QQ_2$, given by $\sqrt d\mapsto\pm\nu$, which
extend to isomorphisms $E_{\Q_j}\z\QQ (i)_{(1+i)}$, given by $\sqrt
d\mapsto\pm\nu$ and $i\mapsto i$. Here $\nu\in\QQ_2$, $\nu\ev 1\pmod
8$ and the $+$ and $-$ signs correspond to $j=1$ and $2$,
respectively. Similarly we define $\nu',\q'_j,\Q'_j$ and we have
isomorphisms $F'_{\q'_j}\z\QQ_2$ and $E'_{\Q'_j}\z\QQ
(i)_{(1+i)}$. These will produce isomorphisms $\psi_j:F'_{\q'_j}\z
F_{\q_j}$ and $\psi'_j:E'_{\Q'_j}\z E_{\Q_j}$. Next we extend
$\Q_j,\Q'_j$ to primes ${\cal Q}_j,{\cal Q}'_j$ of $L,L'$. 

We denote by $\oo_\q$ the ring of integers in $F_\q$ and by $\tilde\q$
the completion of $\q$ in $F_\q$, $\tilde\q =\q\oo_\q$. Similarly for
$\q'$, $\Q$, $\Q'$ and, in the case $d\ev d'\ev 1\pmod 8$, for
$\q_j,\q'_j,\Q_j,\Q'_j$.
\eff

If $d\ev 2,5\pmod 8$ then $\psi'(A'_1)=\psi'((x'-y'\sqrt{d'}i)(u'+v'i)\1
)=(x'-y''\sqrt d)(u'+v'i)\1$, where $y''=\nu y'$. (Recall that $\psi
(\sqrt{d'})=\nu\sqrt d$.) Note that $\nu\ev 1\pmod 8$ so $y'\ev y''\pmod
8$. Also the odd part of $y''$ is $\ev 1\pmod 4$, same as for $y'$,
and also $y''$ is odd, $2\| y''$ or $4\mid y''$ iff $y'$ is so. We
have $x'^2+dy''^2=x'^2+d\nu^2y'^2=x'^2+d'y'^2=p'$.

\blm Under the assumptions from 2.6 and 2.8 we have:

(i) $\frac{u+vi}{u'+v'i}=(1+4i)^{s_1}5^{t_1}B_1^4$, where $B_1\in
1+2\tilde\Q$, $t_1\ev\frac{u'-u}4\pmod 2$ and $s_1\in\{ 0,1\}$,
$s_1\ev\frac{v'-v}4\pmod 2$.
Also $\frac{u'-u}4\ev\frac{p'-p}8\pmod 2$. 

(ii) If $d\ev 5\pmod 8$ and $x,x'$ are odd then $\frac{x'-y''\sqrt di}{x-y\sqrt
di}=(1+4\sqrt di)^{s_2}5^{t_2}B_2^4$, where $B_1\in 1+2\tilde\Q$,
$t_2\ev\frac{x'-x}4\pmod 2$ and $s_2\in\{ 0,1\}$, $s_2\ev\frac{y'-y}4\pmod
2$. 
Also $\frac{x'-x}4\ev\frac{p'-p}8\pmod 2$. 

(iii) If $d\ev 5\pmod 8$ and $x,x'$ are even then $\frac{x'-y''\sqrt di}{x-y\sqrt
di}=5^{t_2}B_2^4$, where $B_1\in 1+2\tilde\Q$, $t_2\ev\frac{y'-y}4\pmod
2$. Also $\frac{y'-y}4\ev\frac{p'-p}8\pmod 2$. 

(iv) If $d\ev 2\pmod 8$ then $\frac{x'-y''\sqrt di}{x-y\sqrt di}=(1+2\sqrt
d)^{2s}5^{t_2}B_2^4$, where $B_1\in 1+2\tilde\Q$,
$t_2\ev\frac{x'-x}4\pmod 2$ and $s\in\{ 0,1\}$, $s\ev\frac{y'-y}4\pmod 2$.
Also $\frac{x'-x}4\ev\frac{p'-p}8\pmod 2$.  
\elm
\pf First note that $v'^2\ev v^2\pmod{16}$. (They are both $\ev
4\pmod{16}$ if $2\| v,v'$ and they are $\ev 0\pmod{16}$ if $4\mid
v,v'$.) Similarly if $d\ev 5\pmod 8$ we have $y^2\ev y''^2\pmod{16}$
if $y,y''$ are even and $x^2\ev x'^2\pmod{16}$ if $x,x'$ are even. The
congruence $y^2\ev y''^2\pmod{16}$ also holds when $d\ev 2\pmod
8$. (See the assumptions of 2.8 and recall that $2\| y''$ or $4\mid
y''$ if $2\| y'$ or $4\mid y'$, respectively.)

We have $p=u^2+v^2$ and $p'=u'^2+v'^2$ so $v'^2\ev v^2\pmod{16}$
implies that $u^2-u'^2\ev p'-p\pmod{16}$. Since $8\mid p'-p$ we have
that $\frac{u'^2-u^2}8$ and $\frac{p'-p}8$ are integers with the same
parities. We also have $u\ev u'\ev 1\pmod 4$ so $2\| u+u'$. Hence
either $2^3\| u'^2-u^2$ and $2^2\| u'-u$ or $16\mid u'^2-u^2$ and
$8\mid u'-u$. In both cases
$\frac{u'-u}4\ev\frac{u'^2-u^2}8\ev\frac{p'-p}8\pmod 2$. 

Similarly $p=x^2+dy^2$ and $p'=x'^2+dy''^2$ so in both cases when $d\ev
5\pmod 8$ and $y,y''$ are even or $d\ev 2\pmod 8$ (so again $y,y''$
are even) the congruence $y^2\ev y''^2\pmod{16}$
implies $x'^2-x^2\ev p'-p\pmod{16}$. Since aslo $x\ev x'\ev 1\pmod 4$
so $2\| x'+x$ we get $\frac{x'-x}4\ev\frac{x'^2-x^2}8\ev\frac{p'-p}8\pmod 2$,
same as above. If $d\ev 5\pmod 8$ and $x,x'$ are even we have
$d(y''^2-y^2)\ev p'-p\pmod{16}$, which implies that
$\frac{p'-p}8,d\frac{y''^2-y^2}8\in\ZZ$ and $\frac{p'-p}8\ev
d\frac{y''^2-y^2}8\ev\frac{y''^2-y^2}8\pmod 2$. Since also $y\ev y''\ev
1\pmod 4$ we have $2\| y''+y$ and so
$\frac{y''-y}4\ev\frac{y''^2-y^2}8\ev\frac{p'-p}8\pmod 2$. But $y'\ev y''\pmod
8$ so $\frac{y'-y}4\ev\frac{y''-y}4\pmod 2$.

So we have proved the second claims of (i)-(iv). We now prove the
first part.  

(i) Suppouse first that $v\ev v'\pmod 8$. We have two cases:

1. $u\ev u'\pmod 8$ or, equivalently, $\frac{u'-u}4$ is even. Then write
$u-u'=8a$, $v-v'=8b$. We have
$(u+vi)-(u'+v'i)=(u-u')+(v-v')i=8(a+bi)$. If $16\mid u-u'+v-v'$,
i.e. if $a+b$ is even, then $a+bi\in\tilde\Q$ so $(u+vi)-(u'+v'i)\in
8\tilde\Q$. Since also $u'+v'i$ is a unit in $\oo_\Q$ we have
$\frac{u+vi}{u'+v'i}-1\in 8\tilde\Q$. By [Lemma 1.6(i)] this implies that
$\frac{u+vi}{u'+v'i}=B_1^4$ for some $B_1\in 1+2\tilde\Q$. Suppose now
that $a+b$ is odd, i.e. that $u-u'+v-v'\ev 8\pmod{16}$. Now $25u'\ev
u'+8\pmod{16}$ and $25v'\ev v'\pmod{16}$. (We have $24\ev 8\pmod{16}$,
$u'$ is odd and $v'$ is even so $24u'\ev 8\pmod{16}$ and $24v'\ev
0\pmod{16}$.) It follows that $25u'\ev u'\ev u\pmod 8$, $25v'\ev v'\ev
v\pmod 8$ and $u-25u'+v-25v'\ev u-u'-8+v-v'\ev 0\pmod{16}$. By the
same reasoning as above we get that $(u+vi)-(25u'+25v'i)\in 8\tilde\Q$
which will imply that $\frac{u+vi}{25u'+25v'i}=B_1^4$, so
$\frac{u+vi}{u'+v'i}=25B_1^4$ for some $B_1\in 1+2\tilde\Q$. Hence,
again, $\frac{u+vi}{u'+v'i}=5^{t_1}B_1^4$ with $t_1$ even. 

2. $u\ev u'+4\pmod 8$, or equivalently, $\frac{u'-u}4$ is odd. Since $u'$
is odd and $v'$ is even we have $4u'\ev 4\pmod 8$ and $4v'\ev 0\pmod
8$ so $5u'\ev u'+4\ev u\pmod 8$ and $5v'\ev v'\ev v\pmod 8$. If we
apply the same reasoning from the case 1. with $u'$ and $v'$ replaced
by $5u'$ and $5v'$ we get $\frac{u+vi}{5u'+5v'i}=5^{t'_1}B_1^4$ with
$t'_1$ even and $B_1\in 1+2\tilde\Q$. It follows that
$\frac{u+vi}{u'+v'i}=5^{t_1}B_1^4$ with $t_1=t'_1+1$ odd, as claimed. 

Suppose now that $v\not\ev v'\pmod 8$. This implies $4\mid v,v'$ and
$v'\ev v+4\pmod 8$, since otherwise $2\| v,v'$ and the odd parts of
$v,v'$ are $\ev 1\pmod 4$ so $v\ev v'\ev 2\pmod 8$ (see 2.8). We have
$(u'+v'i)(1+4i)=u_1+v_1i$, where $u_1=u'-4v'$ and $v_1=v'+4u'$. Since
$v'$ is even and $u'$ odd we have $u_1\ev u'\pmod 8$ and $v_1\ev
v'+4\ev v\pmod 8$. By applying the reasoning from the case $v\ev
v'\pmod 8$ with $u',v'$ replaced by $u_1,v_1$ we get
$\frac{u+vi}{u_1+v_1i}=5^{t_1}B_1^4$, where
$t_1\ev\frac{u_1-u}4\ev\frac{u'-u}4\pmod 2$ and $B_1\in 1+2\tilde\Q$. But
this implies $\frac{u+vi}{u'+v'i}=(1+4i)5^{t_1}B_1^4$. 

(ii) The proof is similar to the proof of (i). 

First note that if $x'-x=8a$ and $y''-y=8b$ and $a+b$ is even then
$a-b\sqrt di\in\tilde\Q$ so $(x'-y''\sqrt di)-(x-y\sqrt di)=8(a-b\sqrt di)\in
8\tilde\Q$. Since also $x-y\sqrt di$ is a unit in $\oo_\Q$ we have
$\frac{x'-y''\sqrt di}{x-y\sqrt di}-1\in 8\tilde\Q$, which implies that
$\frac{x'-y''\sqrt di}{x-y\sqrt di}=B_2^4$, with $B_2\in 1+2\tilde\Q$. From
here the proof in the case $y\ev y''\pmod 8$
follows same as in the case $v\ev v'\pmod 8$ of (i). 

Suppose now that $y\not\ev y''\pmod 8$. This implies that $4\mid
y,y''$ and $y\ev y''+4\pmod 8$. Note that $(x-y\sqrt di)(1+4\sqrt
di)=x_1-y_1\sqrt d$, where $x_1=x+4dy$ and $y_1=y-4x$. Since $x$ is odd
and $y$ is even we have $x_1\ev x'\pmod 8$ and $y_1\ev y+4\ev y''\pmod
8$. By the same reasoning as in the case $y\ev y''\pmod 8$, with
$x_1,y_1$ replacing $x,y$, we get $\frac{x'-y''\sqrt di}{x_1-y_1\sqrt
di}=5^{t_2}B_2^4$, where $B_2\in 1+2\tilde\Q$ and
$t_2\ev\frac{x'-x_1}4\ev\frac{x'-x}4\pmod 2$. It follows that $\frac{x'-y''\sqrt
di}{x-y\sqrt di}=(1+4\sqrt di)5^{t_2}B_2^4$, as claimed. 

(iii) First note that we always have $x\ev x'\pmod 8$. Indeed, if $2\|
x,x'$ then $x\ev x'\ev 2\pmod 8$, while if $4\mid x,x'$ then $\frac
x4\ev\frac{x'}4\pmod 2$ by 2.8(1). If also $y\ev y''\pmod 8$ and we have
$x'-x=8a$, $y'-y=8b$ with $a+b$ even then, same as in the proof of
(ii), $\frac{x'-y''\sqrt di}{x-y\sqrt di}=B_2^4$,  with $B_2\in
1+2\tilde\Q$. From here the proof follows same as in the case $v\ev
v'\pmod 8$ of (i), with the roles of $u,v,u',v'$ being played by
$y'',x',y,x$, respectively.

(iv) Assume first that $y\ev y'\pmod 8$ so $y\ev y''\pmod 8$. We have
two cases:

1. $x\ev x'\pmod 8$ or, equivalently, $\frac{x'-x}4$ is even. If $x\ev
x'\pmod{16}$ then $x'-x\in 16\oo_\Q\sb 8\tilde\Q$ and, since $8\mid
y''-y$ and $\sqrt di\in\tilde\Q$, we also have $(y''-y)\sqrt di\in
8\tilde\Q$. It follows that $(x'-y''\sqrt di)-(x-y\sqrt di)\in 8\tilde\Q$,
which, by the same reasoning as for (ii), will imply that $\frac{x'-y''\sqrt
di}{x-y\sqrt di}=B_2^4$ with $B_2\in 1+2\tilde\Q$. Suppose now that
$x\ev x'+8\pmod{16}$. Since $x$ is odd and $y$ is even we have $24x\ev
8\pmod{16}$ and $24y\ev 0\pmod{16}$ so $25x\ev x+8\ev x'\pmod{16}$ and
$25y\ev y\ev y''\pmod{16}$. By the same reasoning as above
$\frac{x'-y''\sqrt di}{25x-25y\sqrt di}=B_ 2^4$ with $B_2\in 1+2\tilde\Q$ and
so $\frac{x'-y''\sqrt di}{x-y\sqrt di}=25B_2^4$. So in both cases $\frac{x'-y''\sqrt
di}{x-y\sqrt di}=5^{t_2}B_2^4$ with $t_2$ even. 

2. $x\ev x'+4\pmod 8$ or, equivalently, $\frac{x'-x}4$ is odd. Our
statement will follows from the case 1. by the same reasoning as in
the proof of (i) in the case when $v\ev v'\pmod 8$. 

Suppose now that $y\not\ev y''\pmod 8$. It follows that $4\mid y,y''$
and $y\ev y''+4\pmod 8$. By the same reasoning as in the case of (ii)
for $y\not\ev y''\pmod 8$ we get $\frac{x'-y''\sqrt di}{x-y\sqrt
di}=(1+4\sqrt di)5^{t_2}{B'_2}^4$ with $B'_1\in 1+2\tilde\Q$ and
$t_2\ev\frac{x'-x}4\pmod 2$. 

We want to prove that $1+4\sqrt di\ev (1+2\sqrt d)^2=1+4\sqrt
d+4d\pmod{8\tilde\Q}$, i.e. that $\sqrt di\ev\sqrt d+d\pmod{2\tilde\Q}$. 

First note that $E_\Q\ap\QQ_2(\sqrt{-1},\sqrt d)$ is a totally ramified
extension of degree $4$ of $\QQ_2$. Thus $2\oo_\Q =\tilde\Q^4$ and the
inertia degree $[\oo_\Q/\tilde\Q :\ZZ/2\ZZ ]$ is $1$. Thus
$\oo_\Q/\tilde\Q =\{\hat 0,\hat 1\}$. We have $\legendre{i-1}{\sqrt
d}^2=\frac{-2i}d$, which is a unit in $\oo_\Q$, so $\frac{i-1}{\sqrt d}$ is a
unit as well. Since $\frac{i-1}{\sqrt d}\notin\tilde\Q$ its image in
$\oo_\Q/\tilde\Q$ is not $\hat 0$ so it must be $\hat 1$. Thus
$\frac{i-1}{\sqrt d}\ev 1\pmod\Q$. Since also $d\in 2\oo_\Q$ we get $\sqrt
di-\sqrt d=d\cdot\frac{i-1}{\sqrt d}\ev d\pmod{2\tilde\Q}$ and so $\sqrt di\ev\sqrt
d+d\pmod{2\tilde\Q}$, as claimed. 

Since $1+4\sqrt di\ev (1+2\sqrt d)^2\pmod{8\tilde\Q}$ and $1+2\sqrt d$ is a
unit in $\oo_\Q$ we have $\frac{1+4\sqrt di}{(1+2\sqrt d)^2}\ev
1\pmod{8\tilde\Q}$ and so $\frac{1+4\sqrt di}{(1+2\sqrt d)^2}={B''_2}^4$ for some
$B''_2\in 1+2\tilde\Q$. It follows that $\frac{x'-y''\sqrt di}{x-y\sqrt
di}=(1+4\sqrt di)5^{t_2}{B'_2}^4=(1+2\sqrt d)^25^{t_2}B_2^4$, where
$B_2=B'_2B''_2$. Since $B'_2,B''_2\in 1+2\tilde\Q$ we have $B_2\in
1+2\tilde\Q$. \qed 

\blm Let $\alpha,\beta,\alpha',\beta'$ be 2-adic integers
s.t. $\alpha^2+\beta^2=\alpha'^2+\beta'^2=1$. Suppose that $\alpha,\alpha'$ are odd,
$\alpha\ev\alpha'\pmod 4$ and $\beta\ev\beta'\pmod 8$. Then in $\QQ (i)_{(1+i)}$ we
have $\frac{a'+\beta'i}{\alpha +\beta i}=5^{2t}B^4$, where $B\ev 1\pmod{2+2i}$ and
$t\in\ZZ$. Same happens if $\beta,\beta'$ are odd, $\beta\ev\beta'\pmod 4$ and
$\alpha\ev\alpha'\pmod 8$.  
\elm
\pf Let $\oo$ be the ring of integers of $\QQ (i)_{(1+i)}$ and let
$\mm$ be its maximal ideal. We have $\mm =(1+i)\oo$ so the condition
$B\ev 1\pmod{2+2i}$ can be written as $B\ev 1\pmod{2\mm}$ or $B\in
1+2\mm$. 

Suppose first that $\alpha,\alpha'$ are odd, $\alpha\ev\alpha'\pmod 4$ and
$\beta\ev\beta'\pmod 8$. We have that $\beta,\beta'$ are even which, together with
$\beta\ev\beta'\pmod 4$, implies that $\beta^2\ev\beta'^2\pmod{16}$. (We have
$\beta^2\ev\beta'^2\ev 4\pmod{16}$ if $2\|\beta,\beta'$ and
$\beta^2\ev\beta'^2\ev 0\pmod{16}$ if $4\mid\beta,\beta'$.) It follows that
$16\mid\beta^2-\beta'^2=\alpha'^2-\alpha^2$. But also $\alpha,\alpha'$ are odd and
$\alpha\ev\alpha'\pmod 4$ so $2\|\alpha'+\alpha$. It follows that $8\mid\alpha'-\alpha$. Let
$\alpha'-\alpha =8a$ and $\beta'-\beta =8b$. Then $(\alpha'+\beta'i)-(\alpha +\beta
i)=8(a+bi)$. If $16\mid\alpha'-\alpha +\beta'-\beta =8(a+b)$, i.e. if $a+b$ is even,
then $a+bi\in\mm$ so $(\alpha'+\beta'i)-(\alpha +\beta i)\in 8\mm$. But $\alpha +\beta i$
is a unit in $\oo$ so we have $\frac{\alpha'+\beta'i}{\alpha +\beta i}-1\in 8\mm$. By
[B, Lemma 1.6(i)] we get $\frac{\alpha'+\beta'i}{\alpha +\beta i}=B^4$ with $B\in
1+2\mm$. If $\alpha'-\alpha +\beta'-\beta\ev 8\pmod{16}$ then, by the same reasoning
from the proof of Lemma 2.10(i), we have $8\mid\alpha'-25\alpha$,
$8\mid\beta'-25\beta$ and $16\mid\alpha'-25\alpha +\beta'-25\beta$ so this time we get
$\frac{\alpha'+\beta'i}{25\alpha +25\beta i}=B^4$  with $B\in 1+2\mm$ so
$\frac{\alpha'+\beta'i}{\alpha +\beta i}=5^2B^4$.

If $\beta,\beta'$ are odd then note that $\frac{\alpha'+\beta'i}{\alpha +\beta
i}=\frac{\beta'-\alpha'i}{\beta -\alpha i}$ so we repeat the reasoning above with
$\alpha,\alpha',\beta,\beta'$ replaced by $\beta,\beta',-\alpha,-\alpha'$. \qed

\blm (i) If $d\ev d'\ev 5\pmod 8$ then $\psi'(A'_1)/A_1=(3+2\sqrt
d)^{2s}5^{2t}B^4$, where $s\in\{ 0,1\}$, $t\in\ZZ$ and $B\in
1+2\tilde\Q$.  

(ii) If $d\ev d'\ev 2\pmod 8$ then $\psi'(A'_1)/A_1=(1+2\sqrt
d)^{2s}5^{2t}B^4$, where $s\in\{ 0,1\}$, $t\in\ZZ$ and $B\in
1+2\tilde\Q$. Moreover $s\ev\frac{y'-y}4\pmod 2$.

(iii) If $d\ev d'\ev 1\pmod 8$ and $j\in\{ 1,2\}$ then
$\psi'_j(A'_1)/A_1=5^{2t}B^4$, where $t\in\ZZ$ and $B\in
1+2\tilde\Q_j$. 
\elm
\pf For (i) and (ii) we have $\psi'(A_1)/A_1=\frac{x'-y''\sqrt di}{x-y\sqrt
di}\cdot\frac{u+vi}{u'+v'i}$. By Lemma 2.10 we get
$\psi'(A_1)/A_1=(1+4i)^{s_1}(1+4\sqrt di)^{s_2}5^{t_1+t_2}B_1^4B_2^4$ if
$d\ev 5\pmod 8$ and $\psi'(A_1)/A_1=(1+4i)^{s_1}(1+2\sqrt
d)^{2s}5^{t_1+t_2}B_1^4B_2^4$ if $d\ev 2\pmod 8$. (If $d\ev 5\pmod 8$
and $x,x'$ are even we put $s_2=0$, which agrees with Lemma
2.10(iii).) Note that $t_1\ev t_2\ev\frac{p'-p}8\pmod 2$ so $t_1+t_2=2t$
with $t\in\ZZ$ We now prove (i) and (ii). 

For (i) we consider the cases a - d of 2.8(1).

a. $2\| x,x'$ and $\frac v4\ev\frac{v'}4\pmod 2$. Since $x,x'$ are even we
have $s_2=0$ and since $v\ev v'\pmod 8$ we have $s_1=0$. 

b. $2\|y,y'$. By 1.2 this implies that $p\ev p'\ev 5\pmod 8$ so $2\|
v,v'$. We have $y\ev y'\ev 2\pmod 8$ so $s_2=0$ and $v\ev v'\ev 2\pmod
8$ so $s_1=0$. 

c. $4\mid y,y'$ and $\frac{y+v}4\ev\frac{y'+v'}4\pmod 2$. By 1.2 we have
$p\ev p'\ev 1\pmod 8$ so $4\mid v,v'$. The congruence
$\frac{y+v}4\ev\frac{y'+v'}4\pmod 2$ implies that either $\frac
y4\ev\frac{y'}4\pmod 2$ and $\frac v4\ev\frac{v'}4\pmod 2$ or $\frac
y4\ev\frac{y'}4+1\pmod 2$ and $\frac v4\ev\frac{v'}4+1\pmod 2$. In the first
case we get $v\ev v'\pmod 8$ and $y\ev y'\pmod 8$ so $s_1=s_2=0$,
while in the second case $v\ev v'+4\pmod 8$ and $y\ev y'+4\pmod 8$ so
$s_1=s_2=1$. 

d. $4\mid x,x'$ and $\frac x4\ev\frac{x'}4\pmod 2$. Since $x,x'$ are even we
have $s_2=0$. By 1.2 we also have $p\ev p'\ev 5\pmod 8$ so $2\| v,v'$,
which implies $v\ev v'\ev 2\pmod 8$ so $s_1=0$. 

So we proved that $s_1=s_2=0$ or $s_1=s_2=1$. In the first case we get
$\psi'(A_1)/A_1=5^{2t}B^4$, where $B=B_1B_2$. Since $B_1,B_2\in
1+2\tilde\Q$ we have $B\in 1+2\tilde\Q$. In the second case
$\psi'(A_1)/A_1=(1+4i)(1+4\sqrt di)5^{2t}B_1^4B_2^4$. 

We note that for $\alpha,\beta\in\oo_\Q$ we have $(1+4\alpha )(1+4\beta )\ev 1+4(\alpha
+\beta )\pmod{16}$ and $(1+4\alpha )^2\ev 1+8\alpha\pmod{16}$. But $16\in
8\tilde\Q$ so these congruences also hold $\mod 8\tilde\Q$. Since
$d\ev 1\pmod 4$ we have $\frac{1+\sqrt d}2\in\oo_{\QQ (\sqrt
d)}\sb\oo_\Q$. Also note that  $(i-1)^2=-2i\in\tilde\Q$ so
$i-1\in\tilde\Q$ so $i\ev 1\pmod{\tilde\Q}$. In conclusion,
$(1+4i)(1+4\sqrt di)\ev 1+4(1+\sqrt d)i=1+8\frac{1+\sqrt d}2i\ev 1+8\frac{1+\sqrt d}2\ev
(1+4\frac{1+\sqrt d}2)^2=(3+2\sqrt d)^2\pmod{8\tilde\Q}$. Since also $3+2\sqrt d$
is a unit in $\oo_\Q$ we have $\frac{(1+4i)(1+4\sqrt di)}{(3+2\sqrt d)^2}\ev
1\pmod{8\tilde\Q}$ and so $\frac{(1+4i)(1+4\sqrt di)}{(3+2\sqrt d)^2}=B_3^4$
with $B_3\in 1+2\tilde\Q$. Hence $\psi'(A_1)/A_1=(1+4i)(1+4\sqrt
di)5^{2t}B_1^4B_2^4=(3+2\sqrt d)^25^{2s}B^4$, where $B=B_1B_2B_3$. Since
$B_1,B_2,B_3\in 1+2\tilde\Q$ we have $B\in 1+2\tilde\Q$. 

(ii) As seen from 1.2 we have $4\mid v,v'$. But in both cases a and b
of 2.8(2) we have $\frac v4\ev\frac{v'}4\pmod 2$ so $v\ev v'\pmod 8$. Hence
$s_1=0$ and we have $\psi'(A_1)/A_1=(1+2\sqrt d)^{2s}5^{2t}B^4$,
where $B=B_1B_2$. But $B_1,B_2\in 1+2\tilde\Q$ so $B\in
1+2\tilde\Q$. 

(iii) We have $\psi'_j(i)=i$ and $\psi'_j(\sqrt{d'})=\pm\nu'$. If 
$\alpha_{\Q_j}$ is the image in $E_{\Q_j}\ap\QQ (i)_{(1+i)}$ of an element
$\alpha\in E$ then $i_{\Q_j}=i$ and $\sqrt d_{\Q_j}=\pm\nu$. Hence
${A_1}_{\Q_j}=\frac{x-y\xi_ji}{u+vi}$ and
$\psi'_j(A_1)=\frac{x'-y'\xi_j'i}{u'+v'i}$, where $\xi_j=\pm\nu$ and
$\xi'_j=\pm\nu'$. Since $\nu\ev\nu'\ev 1\pmod 8$ we have
$\xi_j\ev\xi'_j\ev\pm 1\pmod 8$. In all cases above the sign $+$
corresponds to $j=1$ and $-$ to $j=2$.

We have $\psi'_j(A_1)/A_1=\frac{\alpha'+\beta'i}{\alpha+\beta i}$, where $\alpha +\beta
i=\frac{x-y\xi_ji}{u+vi}$ and $\alpha'+\beta'i=\frac{x'-y'\xi'_ji}{u'+v'i}$. We
have $x^2+y^2\xi^2=x^2+dy^2=p$ and $u^2+v^2=p$. It follows that $\alpha
=(xu-yv\xi_j)p\1$, $\beta =(-xv-yu\xi_j)p\1$ and $\alpha^2+\beta^2=1$. Similarly
for $\alpha',\beta'$. We now prove that $\alpha,\beta,\alpha',\beta'$ satisfy the
conditions of Lemma 2.11. We already have $\alpha^2+\beta^2=\alpha'^2+\beta'^2=1$
and we still have to prove the congruences from Lemma 2.11. From 1.2
we have that if $2\| xy$ and $2\| x'y'$ then $p\ev p'\ev 5\pmod 8$ so
$2\| v$ and $2\| v'$. If $4\mid xy$ and $4\mid x'y'$ then $p\ev 
p'\ev 1\pmod 8$ so $4\mid v$ and $4\mid v'$. In both cases $p\ev
p'\pmod 8$ so $p\1\ev {p'}\1\pmod 8$. It also follows that $x,y,u,v$
are $\ev\mod 4$ to $x',y',u',v'$. (See also 2.8 and the hypothesis.)
Since $\xi_j\ev\xi'_j\pmod 4$ and $p\1\ev {p'}\1\pmod 4$ we have
$\alpha\ev\alpha'\pmod 4$ and $\beta\ev\beta'\pmod 4$. So we still have to prove
that either $\alpha,\alpha'$ are even and $\alpha\ev\alpha'\pmod 8$ or $\beta,\beta'$ are
even and $\beta\ev\beta'\pmod 8$. Note that if $x$ is even then
$\alpha=(xu-yv\xi_j)p\1$ is even, while if $y$ is even then $\beta
=(-xv-yu\xi_j)p\1$ is even. Same for $\alpha',\beta'$. We now consider the
cases a-d of 2.8(3):

a. $2\| x,x'$. It follows that $2\| v,v'$. By hypothesis we have $\frac
x2\ev\frac{x'}2\ev\frac v2\ev\frac{v'}2\ev 1\pmod 4$. Together with $y\ev y'\ev
u\ev u'\ev 1\pmod 4$, $\xi_j\ev\xi'_j\pmod 4$ and $p\1\ev {p'}\1\pmod
4$, this implies that $\frac\alpha 2=(\frac x2u-y\frac v2\xi_j)\p\1\ev
(\frac{x'}2u'-y'\frac{v'}2\xi'_j){\p'}\1 =\frac{\alpha'}2\pmod 4$ so $\alpha\ev\alpha'\pmod
8$. 

b. $2\| y,y'$. Similarly as in case a. we have $\frac y2\ev\frac{y'}2\ev\frac
v2\ev\frac{v'}2\ev 1\pmod 4$, $x\ev x'\ev u\ev u'\ev 1\pmod 4$,
$\xi_j\ev\xi'_j\pmod 4$ and $p\1\ev{p'}\1\pmod 4$. Hence $\frac\beta 2=(-x\frac
v2-\frac y2u\xi_j)p\1\ev (-x'\frac{v'}2-\frac{y'}2u'\xi'_j){p'}\1 =\frac{\beta'}2\pmod
4$ so $\beta\ev\beta'\pmod 8$. 

c. $4\mid x,x'$ and $\frac{x+v}4\ev\frac{x'+v'}4\pmod 2$. We have $4\mid
v,v'$ and we write $\frac\alpha 4=(\frac x4u+y\frac v4\xi_j)\p\1$. But
$y,u,\xi_j,\p\1$ are odd so $\frac\alpha 4\ev\frac{x+v}4\pmod 2$. Similarly
$\frac{\alpha'}4\ev\frac{x'+v'}4\pmod 2$ so $\frac\alpha 4\ev\frac{\alpha'}4\pmod 2$ so
$\alpha\ev\alpha'\pmod 8$. 

d. $4\mid y,y'$ and $\frac{y+v}4\ev\frac{y'+v'}4\pmod 2$. Same as in case
c. we have $4\mid v,v'$ and since $x,u,\xi_j,\p\1$ are odd we have
$\frac\beta 4=(-x\frac v4-\frac y4u\xi_j)p\1\ev\frac{y+v}4\pmod 2$. Similarly
$\frac{\beta'}4\ev\frac{y'+v'}4\pmod 2$ so $\frac\beta 4\ev\frac{\beta'}4\pmod 2$ so
$\beta\ev\beta'\pmod 8$.

By Lemma 2.11 we have $\psi'_j(A_1)/A_1=\frac{\alpha'+\beta'i}{\alpha +\beta
i}=5^{2t}B^2$ with $t\in\ZZ$ and $B\ev 1\pmod{2+2i}$. But $1+i$
generates the ideal $\tilde\Q_j$ so
$(2+2i)\oo_{\Q_j}=2\tilde\Q_j$. Thus $B\ev 1\pmod{2\tilde\Q_j}$. \qed

\blm (i) If $d\ev 5\pmod 8$ then $\legendre{\varepsilon,5}\q =\legendre{\varepsilon,3+2\sqrt d}\q
=1$. 

(ii) If $d\ev 2\pmod 8$ then $\legendre{\varepsilon,5}\q =1$ and $\legendre{\varepsilon,1+2\sqrt d}\q
=-1$. 
\elm
\pf We have $5\in\QQ_2$ and in both cases $d\ev 5$ or $2\pmod 8$ we
have $\N_{F_\q/\QQ_2}(\varepsilon )=-1$. It follows that $\legendre{\varepsilon,5}\q
=\legendre{-1,5}2=1$ so we have proved the first part of (i) and (ii). We
now prove the second part.

(i) We have that $\frac{1+\sqrt d}2$ is an algebraic integer so $3+2\sqrt
d=1+4\frac{1+\sqrt d}2\ev 1\pmod 4$. It follows that $3+\sqrt d$ is either a
square or a unit with quadratic defect $4\oo_\q$. (See [OM, \S 63A].)
Since also $\varepsilon$ is a unit in $F_\q$ we have $\legendre{\varepsilon,3+2\sqrt d}\q
=1$. (See [OM, 63:11a].) 

(ii) First recall some properties of the Hilbert symbol. If $\alpha\in
F_\q^\ti$ then $\legendre{\alpha,-\alpha}\q =1$ so $\legendre{\alpha,\alpha}\q
=\legendre{\alpha,-1}\q$. If moreover $\alpha\neq 1$ then $\legendre{\alpha,1-\alpha}\q =1$.

If $X,Y\in F_\q$ s.t. $XY(X+Y)\neq 0$ then $1-\frac X{X+Y}=\frac Y{X+Y}$ so
$1=\legendre{\frac X{X+Y},\frac Y{X+Y}}\q
=\legendre{X,Y}\q\legendre{X,X+Y}\q\legendre{X+Y,Y}\q\legendre{X+Y,X+Y}\q=
\legendre{X,Y}\q\legendre{X+Y,-XY}\q$. (We have
$\legendre{X+Y,X+Y}\q=\legendre{X+Y,-1}\q$.) Hence $\legendre{X,Y}\q
=\legendre{X+Y,-XY}\q$. 

Let now $\alpha,\beta\in F_\q$, $\alpha\neq -1,0$, $\beta\neq -1$, $\alpha\beta\neq 1$. We
have $\legendre{-\alpha,1+\alpha}\q =1$ so $\legendre{1+\alpha,1+\beta}\q =\legendre{1+\alpha,-\alpha(1+\beta
)}\q$. We now use the property above for $X=1+\alpha$ and $Y=-\alpha(1+\beta
)=-\alpha -\alpha\beta$. We get $\legendre{1+\alpha,1+\beta}\q =\legendre{1-\alpha\beta,-\alpha (1+\alpha )(1+\beta
)}\q$.  

We use now this property for $\alpha,\beta$ given by $1+\alpha =\varepsilon$ and $1+\beta
=1+2\sqrt d$. We claim that $1-\alpha\beta$ is a unit of quadratic defect
$4\oo_\q$ and $-\alpha (1+\alpha )(1+\beta )$ is a prime in $\oo_\q$. By [OM,
63:11a] this will imply $-1=\legendre{1-\alpha\beta,-\alpha (1+\alpha)(1+\beta )}\q
=\legendre{1+\alpha,1+\beta}\q =\legendre{\varepsilon,1+2\sqrt d}\q$, as claimed. 

First note that $F_\q$ is a ramified extension of degree $2$ of
$\QQ_2$ so $2\oo_\q=\tilde\q^2$ and $\sqrt d$ is a prime of
$\oo_\q$. Also the inertia index $[\oo_\q/\tilde\q:\ZZ/2\ZZ ]$ is
$1$.  

Since $d\ev 2\pmod 8$ we are in the case of Conjecture 3, when $b\ev
2\pmod 4$. We have $\alpha =\frac b2-1+\sqrt d$ and $\beta =2\sqrt d$. Now $\sqrt d$ is a
prime and $2\mid\frac b2-1$ so $\frac b2-1\in\tilde\q^2$ so $\alpha =\frac b2-1+\sqrt
d$ is a prime in $\oo_\q$. Since also $1+\alpha =\varepsilon$ and $1+\beta =1+2\sqrt d$
are units in $\oo_\q$, $-\alpha (1+\alpha )(1+\beta )$ will be a prime as well. 

We have $\alpha\beta =2\sqrt d(\frac b2 -1+\sqrt d)=(b-2)\sqrt d +2d=4\rho$, where $\rho
=\frac d2+\frac{b-2}4\sqrt d$. But $\frac d2\ev 1\pmod 2$ so $\frac d2\ev
1\pmod{\tilde\q}$ and we also have $\frac{b-2}4\sqrt d\in\tilde\q$ so
$\rho\ev 1\pmod{\tilde\q}$. 

We have $1-\alpha\beta =1-4\rho$ so $1-\alpha\beta$ can only be either a square or a
unit of quadratic defect $4\oo_\q$. Assume that it is a square so
$1-4\rho =\gamma^2$ for some $\gamma\in F_\q$. Since $4\mid\gamma^2-1$ we have
$2\mid\gamma -1$ or $2\mid\gamma +1$. But if $2\mid\gamma +1$ then also $2\mid\gamma
+1-2=\gamma -1$. So always $2\mid\gamma -1$. Then $\gamma =1+2\delta$ with
$\delta\in\oo_\q$. From $1-4\rho =(1+2\delta )^2$ we get $\delta^2+\delta +\rho
=0$. But $\rho\ev 1\pmod{\tilde\q}$ so by considering classes
$\mod\tilde\q$ we get $\hat\delta^2+\hat\delta +\hat 1=\hat 0$. But this is
impossible since $X^2+X+\hat 1$ has no roots in
$\oo_\q/\tilde\q\ap\ZZ/2\ZZ$. Hence $1-\alpha\beta$ must be a unit of
quadratic defect $4\oo_\q$, as claimed. \qed

\blm If $d,d'$ are odd then $\eta =\eta'$. If $d,d'$ are even then
$\eta =(-1)^{\frac{y-y'}4}\eta'$. 
\elm
\pf We first treat the case when $d\ev d'\ev 5$ or $2\pmod 8$. We have
$\eta=\chi (\phi )$, where $\phi =\legendre{\varepsilon,L/F}\q$. If $\cal Q$ is a
prime of $L$ staying over $\Q$, so over $\q$, then $\phi =(\varepsilon,L_{\cal
Q}/F_\q )$. We have $\phi (i)=\legendre{\varepsilon,-1}\q i$. But $-1\in\QQ_2$ and
$\N_{F_\q/\QQ_2}\varepsilon =-1$ so $\legendre{\varepsilon,-1}\q=\legendre{-1,-1}2=-1$ so $\phi
(i)=-i$. It follows that $\phi_{|E}={\tau_3}_{|E}$ so
$\phi\in\tau_3Gal(L/E)=\tau_3\la\s\ra$ so $\phi =\tau_3\s^k$ for some
$k\in\ZZ_4$. It follows that $\eta =\chi (\phi )=i^k$. We have
$\tau_3(\alpha_1)=\alpha_3=\alpha_1\1$ so $\phi
(\alpha_1)=\tau_3\s^k(\alpha_1)=\tau_3(i^k\alpha_1)=i^{-k}\alpha_1\1$.

We repeat the reasoning for $b',p'$. Denote by $\phi',{\cal Q}',k'$
the $\phi,{\cal Q},k$ corresponding to $b',p'$. We have to prove that
$i^k=i^{k'}$ or $i^k=(-1)^{\frac{y-y'}4}i^{k'}$, corresponding to $d\ev
5$ or $2\pmod 8$, respectively. 

We have the isomorphism $\psi :F'_{\q'}\z F_\q$ given by
$\sqrt{d'}\mapsto\nu\sqrt d$, which extends to an isomorphism $\psi'
:E'_{\Q'}\z E_\Q$ given by $\sqrt{d'}\mapsto\nu\sqrt d$, $i\mapsto i$. By
Lemma 2.12 we have $\psi'(A_1')/A_1=C^2B^4$, where $B\in 1+2\tilde\Q$
and $C=(3+2\sqrt d)^s5^t$ or $(1+2\sqrt d)^s5^t$, corresponding to $d\ev 5$
or $2\pmod 8$, respectively. In both cases $C\in\QQ (\sqrt d)=F$. We
consider the extension $L_{\cal Q}(\lambda )$ of $L_{\cal Q}$ given by
$\lambda^2=C$. Since both $L_{\cal Q}$ and $F_\q (\lambda )\ap F_\q (\sqrt C)$ are
abelian extensions of $F_\q$, so is $L_{\cal Q}(\lambda )$. 

Now $L'_{\cal Q'}=E'_{\Q'}(\alpha_1')$ and $\alpha_1'$ is a root of
$X^4-A_1'\in E'_{\Q'}[X]$. The roots of $\psi'
(X^4-A_1')=X^4-C^2B^4A_1\in E_\Q [X]$ are $i^l\lambda B\alpha_1$ with $0\leq
l\leq 3$ and they all belong $L_{\cal Q}(\lambda )$. Therefore the
isomorphism $\psi':E'_{\Q'}\z E_{\Q}$ extends to a morphism
$\psi'':L'_{\cal Q'}\z L_{\cal Q}(\lambda )$ given by $\alpha'_1\mapsto i^l\lambda
B\alpha_1$ for some $l$. We have:
$$\psi''((\varepsilon',L'_{\cal Q'}/F'_{\q'})(\alpha'_1))=(\psi (\varepsilon'),L_{\cal
Q}(\lambda )/F_\q )(\psi''(\alpha'_1)).$$

The left hand side can be written as
$\psi''(\phi'(\alpha'_1))=\psi''(i^{-k'}{\alpha'_1}\1 )=i^{-k'}(i^l\lambda
B\alpha_1)\1=i^{-k'-l}\lambda\1 B\1\alpha_1\1$.

For the right hand side we have $(\psi (\varepsilon'),L_{\cal Q}(\lambda
)/F_\q )=(A^4\varepsilon,L_{\cal Q}(\lambda )/F_\q )=(\varepsilon,L_{\cal Q}(\lambda )/F_\q
)$. (If $\theta :=(A,L_{\cal Q}(\lambda )/F_\q )$ then
$\theta_{|L_{\cal Q}}\in Gal(L_{\cal Q}/F_\q )\sbq
Gal(L/F)\ap\ZZ_2\ti\ZZ_4$ and $\theta_{|F_\q (\lambda )}\in Gal(F_\q (\lambda
)/F_\q )\hookrightarrow\ZZ_2$. It follows that $\theta_{|L_{\cal
Q}}^4=1_{L_{\cal Q}}$ and $\theta_{|F_\q (\lambda )}^4=1_{F_\q (\lambda )}$ so
$\theta^4=1_{L_{\cal Q}(\lambda )}$.) Since also $\psi''(\alpha'_1)=i^l\lambda
B\alpha_1$, the right hand side equals $(\varepsilon,L_{\cal Q}(\lambda )/F_\q )(i^l\lambda
B\alpha_1)=i^{-k-l}\legendre{\varepsilon,C}\q\lambda\phi(B)\alpha_1\1$. (We have
$\lambda^2=C\in F_\q$ so $(\varepsilon,L_{\cal Q}(\lambda )/F_\q )(\lambda )=\legendre{\varepsilon,C}\q\lambda$
and $(\varepsilon,L_{\cal Q}(\lambda )/F_\q )(i^lB\alpha_1)=\phi (i^lB\alpha_1)=i^{-l}\phi
(B)i^{-k}\alpha_1\1$. 

In conclusion, $i^{-k'-l}\lambda\1 B\1\alpha_1\1 =(-1)^{-k-l}\legendre{\varepsilon,C}\q\lambda\phi
(B)\alpha_1\1$, so $i^{k-k'}\legendre{\varepsilon,C}\q =\lambda^2B\phi (B)=CB\phi (B)$. But
$B$ belongs to $1+2\tilde\Q$ and so does its conjugate $\phi
(B)$. Also $C\in 1+2\tilde\Q$. Indeed, if $d\ev 5\pmod 8$ then $5$ and
$3+2\sqrt d$ belong to $1+4\oo_\Q$ so $C=(3+2\sqrt d)^s5^t\in 1+4\oo_\Q\sb
1+2\tilde\Q$, while if $d\ev 2\pmod 8$ then $5,1+2\sqrt d\in 1+2\tilde\Q$
so $C=(1+2\sqrt d)^s5^t\in 1+2\tilde\Q$. It follows that
$i^{k-k'}\legendre{\varepsilon,C}\q\in 1+2\tilde\Q$. Since
$i^{k-k'}\legendre{\varepsilon,C}\q\in\mu_4$ this implies
$i^{k-k'}\legendre{\varepsilon,C}\q=1$. (If $\xi\in\mu_4$ and $\xi\neq 1$ then $\xi
-1\mid 2$ and so $\xi -1\notin 2\tilde\Q$.) It follows that
$i^k=\legendre{\varepsilon,C}q i^{k'}$. If $d\ev 5\pmod 8$ then $C=(3+2\sqrt d)^s5^t$,
which, by Lemma 2.13(i), implies $\legendre{\varepsilon,C}\q =1$. If $d\ev 5\pmod 8$
then $C=(3+2\sqrt d)^s5^t$, which, by Lemma 2.13(i), implies
$\legendre{\varepsilon,C}\q =(-1)^s$. But $s\ev\frac{y-y'}4\pmod 2$, so
$(-1)^s=(-1)^{\frac{y-y'}4}$. Hence the desired result. 

Suppose now that $d\ev 1\pmod 8$. By Lemma 2.7(ii) we have $\eta
=\eta_1\eta_2$, where
$\eta_j=\chi\left(\legendre{a_j,L/F}{\q_j}\right)$. Here $a_1=1$,
$a_2=-1$ if $16\mid b$ and $a_1=5$, $a_2=3$ if $16\nmid b$. Similarly
for $b',p'$. We denote by $a'_1,a'_2,\eta'_1,\eta'_2$ the
$a_1,a_2,\eta_1,\eta_2$ corresponding to $b',p'$ Recall that by 2.6 we
have either both $16\mid b$ and $16\mid b'$ or both $16\nmid b$ and
$16\nmid b'$ so $a'_j=a_j$. We will prove that $\eta_j=\eta'_j$ for
$j=1,2$. 

We fix now $j\in\{ 1,2\}$. By Lemma 2.12(iii) we have
$\psi'_j(A_1')/A_1=C^2B^4$ for some $B\in 1+2\tilde\Q$, where
$C=5^t$. We proceed similarly as in the cases $d\ev 5,2\pmod 8$. Let
$\phi= \legendre{a_j,L/F}{\q_j}$ and $\phi'=\legendre{a_j,L'/F'}{\q'_j}$. If
$\cal Q,Q'$ are pimes of $L,L'$ staying over $\Q_j,\Q'_j$ then $\phi
=(a_j,L_Q/F_\q )$ and $\phi'=(a_j,L'_{Q'}/F'_{\q'})$. If $L_{\cal
Q}(\lambda )$ is an extension of $L_{\cal Q}$ with $\lambda^2=C$ then the
isomorphism $\psi'_j:E'_{\q'_j}\z E_{\q_j}\ap\QQ (i)_{(1+i)}$ will
extend to a morphism $\psi''_j:L'_{\cal Q}\z L_{\cal Q'}(\lambda )$
satisfying $\psi''_j(\alpha'_1)=i^l\lambda B\alpha_1$ for some $l$. We have 
$$\psi''_j((a_j,L'_{\cal Q'}/F'_{\q'_j})(\alpha'_1))=(a_i,L_{\cal Q}(\lambda
)/F_{\q_j})(\psi''_j(\alpha'_1)).$$

If $j=2$ then $a_2=-1$ or $3$. In both cases
$\legendre{a_2,-1}{\q_2}=\legendre{a_2,-1}2=-1$ so $\phi
(i)=-i$. (Recall that $F_{\q_2}\ap\QQ_2$.) By the same reasoning as in
the previous case we get $\phi =\tau_3\s^k$ and
$\phi'=\tau'_3\s'^{k'}$ for some integers $k,k'$ and we have to prove
that $i^k=i^{k'}$. Same as before we get
$i^{k-k'}\legendre{C,a_2}{\q_2}=CB\phi (B)$. Again $B,\phi (B)\in
1+2\tilde\Q_2$ and $5\in 1+4\oo_{\Q_2}\sb 1+2\tilde\Q_2$ so $C=5^t\in
1+2\tilde\Q_2$. It follows that $i^{k-k'}\legendre{C,a_2}{\q_2}\in
1+2\tilde\Q_2$. Together with $i^{k-k'}\legendre{C,a_2}{\q_2}\in\mu_4$,
this implies $i^{k-k'}\legendre{C,a_2}{\q_2}=1$. But in both cases
$a_2=-1,3$ we have $\legendre{5,a_2}{\q_2}=\legendre{5,a_2}2=1$. Since $C=5^t$
we get $\legendre{C,a_2}{\q_2}=1$ so $i^{k-k'}=1$ and we are done. 

Suppose now that $j=1$. If $a_1=1$ then $\eta_1=\eta'_1=1$ so we may
assume that $a_1=5$. Since $\legendre{5,-1}{\q_1}=\legendre{5,-1}2=1$
we have $\phi (i)=i$ and so $\phi_{|E}=1_E$. Thus $\phi\in
Gal(L/E)=\la\s\ra$. If $\phi =\s^k$ then $\eta_1=\chi (\phi )=i^k$ and
$\phi (\alpha_1)=i^k\alpha_1$. Similarly $\phi'=\s'^{k'}$, $\eta'_2=i^{k'}$
and $\phi'(\alpha'_1)=i^{k'}\alpha'_1$. We now calculate the left and right
hand sides of the equality above. 

The left hand side equals
$\psi''_1(\phi'(\alpha'_1))=\psi''_1(i^{k'}\alpha'_1)=i^{k'+l}\lambda B\alpha_1$.

The right hand side equals $(5,L_{\cal Q}(\lambda )/F_{\q_1})(i^l\lambda
B\alpha_1)=\legendre{5,C}{\q_1}i^{k+l}\lambda\phi (B)\alpha_1$. (We have $\lambda^2=C$ so
$(5,L_{\cal Q}(\lambda )/F_{\q_1})(\lambda )=\legendre{5,C}{\q_1}\lambda$ and $i^lB\alpha_1\in
L_{\cal Q}$ so $(5,L_{\cal Q}(\lambda )/F_{\q_1})(i^lB\alpha_1)=(5,L_{\cal
Q}/F_{\q_1})(i^lB\alpha_1)=\phi (i^lB\alpha_1)=i^l\phi (B)i^k\alpha_1$.) 

So $i^{k'+l}\lambda B\alpha_1=\legendre{5,C}{\q_1}i^{k+l}\lambda\phi (B)\alpha_1$ which implies
$i^{k-k'}\legendre{5,C}{\q_1}=B\phi (B)\1$. Since $B,\phi (B)\in 1+2\tilde\Q_j$ we
have $i^{k-k'}\legendre{5,C}{\q_1}\in 1+2\tilde\Q_j$, which, together with
$i^{k-k'}\legendre{5,C}{\q_1}\in\mu_4$ implies $i^{k-k'}\legendre{5,C}{\q_1}=1$. But
$\legendre{5,5}{\q_1}=\legendre{5,5}2=1$ and $C=5^t$ so $\legendre{5,C}{\q_1}=1$ so
$i^{k-k'}=1$. \qed

\blm With the conventions of 2.6 and 2.8, the one of the Conjectures
1-4 is true for the pair $b,p$ iff it is true for $b',p'$.
\elm
\pf Conjectures 1-4 are true for a pair $b,p$ iff $\eta =\xi$, where
$\xi$ is the value predicted for $\eta$ in Lemma 2.5. Similarly, if we
denote by $\xi'$ the value predicted for $\eta'$, the Conjectures 1-4
are true for $b',p'$ iff $\eta'=\xi'$. In order to prove that the two
statements are equivalent it is enough to show that
$$\frac{\xi'}\xi =\frac{\eta'}\eta =\begin{cases}(-1)^{\frac{y-y'}4}&\text{if
}d,d'\text{ are even}\\1&\text{if }d,d'\text{ are
  odd}\end{cases}\text{~~~(see Lemma 2.14)}.$$

We consider the parts (i)-(iv) of Lemma 2.5, which correspond to
Conjectures 1-4. 

(i) By 2.6(1) we have $b\ev b'\pmod 8$ so Lemma 2.5(i) implies that:
$$\xi'/\xi =\begin{cases}(-1)^{\frac{v'-v}4}&\text{if }2\| x,x'\\
1&\text{ if }2\| y,y'\\ (-1)^{\frac{y'+v'}4-\frac{y+v}4}&\text{if }4\mid
y,y'\\ (-1)^{\frac{x'-x}4}&\text{if }4\mid x,x'.\end{cases}$$
By 2.8(1) in all four cases $\xi'/\xi =1$. 

(ii) By 2.6(3) we have $b\ev b'\pmod{16}$ so $\frac{b'-b}8$ is even. So
by Lemma 2.5(ii) we have
$$\xi'/\xi
=\begin{cases}(-1)^{\frac{b'-b}8+\frac{v'-v}4}=(-1)^{\frac{v'-v}4}&\text{if
}2\| x,x'\\ 1&\text{if }2\| y,y'\\ (-1)^{\frac{y'+v'}4-\frac{y+v}4}&\text{if
}4\mid y,y'\\ (-1)^{\frac{b'-b}8+\frac{x'-x}4}=(-1)^{\frac{x'-x}4}&\text{if
}4\mid x,x'.\end{cases}$$ 
By 2.8(1) in all four cases $\xi'/\xi =1$. 

(iii) Each of the cases $4\nmid xy$ and $4\mid xy$ from Lemma 2.5(iii)
can be split in to subcases, according to the parities of $x,y$. If
$2\| x$ then $y$ is odd so $\xi =(-1)^{(\frac b8-1)y}=(-1)^{\frac b8-1}$. If
$2\| y$ then $\xi =(-1)^{(\frac b8-1)y}=1$. If $4\mid x$ then $y$ is odd
so $xy\ev x\pmod 8$ so $\frac{xy+v}4\ev\frac{x+v}4\pmod 2$. Also $\frac
b8y\ev\frac b8\pmod 2$. Thus $\xi =(-1)^{\frac{xy+v}4+\frac
b8y}=(-1)^{\frac{x+v}4+\frac b8}$. Similarly, if $4\mid y$ we get that
$\frac{xy+v}4\ev\frac{y+v}4\pmod 2$ and, since also $\frac b8y$ is even, we get
$\xi =(-1)^{\frac{xy+v}4+\frac b8y}=(-1)^{\frac{y+v}4}$. Similarly for
$\xi'$. We get $$\xi'/\xi =\begin{cases}(-1)^{\frac{b'-b}8}&\text{if
}2\|x,x'\\ 1&\text{if }2\|
y,y'\\ (-1)^{\frac{x'+v'}4-\frac{x+v}4+\frac{b'-b}8}&\text{if }4\mid
x,x'\\ (-1)^{\frac{y'+v'}4-\frac{y+v}4}&\text{if }4\mid y,y'.\end{cases}$$ 
By 2.6(4) we have $b\ev b'\pmod{16}$ so $\frac{b'-b}8$ is even. Together
with 2.8(3), this implies, in all four cases, that $\xi'/\xi =1$. 

(iv) By 2.6(2) we have $b\ev b'\pmod{32}$ so $\frac{b'-b}4$ is even. By
2.8(2) we have that $\frac{v'-v}4$ is even. So Lemma 2.5(iv) implies
$$\xi'/\xi =\begin{cases}(-1)^{\frac{b'-b}4+\frac{v'-v}4}=1&\text{if }2\|
y,y'\\ (-1)^{\frac{y'-y}4+\frac{v'-v}4}=(-1)^{\frac{y'-y}4}&\text{if }4\mid
y,y'.\end{cases}$$ 
But if $2\| y,y'$ then $y\ev y'\ev 2\pmod 8$ so $\frac{y'-y}4$ is
even. Thus $\xi'/\xi =(-1)^{\frac{y'-y}4}$ holds in this case as
well. \qed 

\bff Let now $b,p$ be a pair satisfying the hypothesis of one of the
Conjectures 1-4. We consider several cases. 

(I) If $d\ev 5\pmod 8$ then either $b$ is odd or $b\ev 4\pmod 8$. We
have three cases: (1) $b\ev\pm 1\pmod 8$, (2) $b\ev\pm 3\pmod 8$, (3)
$b\ev 4\pmod 8$. If two integers $b,b'\neq 0$ are in the same case
then either $b$ and $b'$ or $-b$ and $b'$ satisfy the conditions of
2.6. For each of these cases we have seven subcases involving the
values of $x,y,u,v$ modulo powers of $2$: 

a. $2\| x$, $\frac v4$ is even.

b. $2\| x$, $\frac v4$ is odd. 

c. $2\| y$

d. $4\mid x$, $\frac x4$ is even.

e. $4\mid x$, $\frac x4$ is odd.

f. $4\mid y$, $\frac{y+v}4$ is even. 

g. $4\mid y$, $\frac{y+v}4$ is odd. 

So we have $3\cdot 7=21$ cases. If the pairs $b.p$ and $b',p'$
satisfying the hypothesis of Conjecture 1 or 2 fall in the same case
then either $b,p$ and $b',p'$ or $-b,p$ and $b',p'$ satisfy the
conditions of both 2.6(1) or (3) and of 2.8(1). 

(II) $d\ev 2\pmod 8$ so $b\ev 2\pmod 4$. We have four cases: (1)
$b\ev\pm 2\pmod{32}$, (2) $b\ev\pm 6\pmod{32}$, (3) $b\ev\pm
10\pmod{32}$ and (4) $b\ev\pm 14\pmod{32}$. 

For each of these cases we have four subcases:

a. $2\| y$, $\frac v4$ is even.

b. $2\| y$, $\frac v4$ is odd. 

c. $4\mid y$, $\frac v4$ is even.

d. $4\mid y$, $\frac v4$ is odd. 

There are $4\cdot 4=16$ cases. If two pairs $b,p$ and $b',p'$
satisfying the hypothesis of Conjecture 3 fall in the same class then
either $b,p$ and $b',p'$ or $-b,p$ and $b',p'$ satisfy the conditions
of both 2.6(2) and 2.8(2). 

(III) $d\ev 1\pmod 8$ so $8\mid b$. We have two cases: (1) $b\ev
0\pmod{16}$ and (2) $b\ev 8\pmod{16}$. 

For each of these cases we have six subcases:

a. $2\| x$.

b. $2\| y$.

c. $4\mid x$, $\frac{x+v}4$ is even.  

d. $4\mid x$, $\frac{x+v}4$ is odd. 

e. $4\mid y$, $\frac{y+v}4$ is even.  

f. $4\mid y$, $\frac{y+v}4$ is odd. 

There are $2\cdot 6=12$ cases. If two pairs $b,p$ and $b',p'$
satisfying the hypothesis of Conjecture 3 fall in the same class then
either $b,p$ and $b',p'$ or $-b,p$ and $b',p'$ satisfy the conditions
of both 2.6(4) and 2.8(3). 
\eff 

{\bf End of the proof } If $b,p$ and $b',p'$ are two pairs satisfying
the hypothesis of one of the Sun's conjectures which fall in one of
the $21+16+12=49$ cases of 2.16 then either $b,p$ and $b',p'$ or
$-b,p$ and $b',p'$ satisfy the conditions of both 2.6 and 2.8. By
Lemma 2.15 the conjecture for the two pairs is equivalent. But by
Lemma 1.3 the conjecture holds for $b,p$ iff it holds for
$-b,p$. Hence in both cases above our conjecture for $b,p$ and for
$b',p'$ is equivalent. Therefore the Conjectures 1-4 need to be
verified only for a set of pairs $b,p$ that cover all $49$ cases of
2.16. We consider the cases (I)-(III) of 2.16.

(I) $b=1,3,4$ cover the cases (1)-(3) of (I). For each of these we
display seven values of $p$ covering the cases a-g from 2.16(I). In
for these primes we also show how they can be written as
$p=x^2+dy^2=u^2+v^2$.  

If $b=1$ then $d=5$ and we take $p=281,41,29,109,61,521,89$. We have
$241=(-6)^2+5\cdot (-7)^2=5^2+16^2$, $41=(-6)+5\cdot 1^2=5^2+4^2$,
$29=(-3)^2+5\cdot 2^2=5^2+2^2$, $109=8^2++5\dot (-3)^2=(-3)^2+10^2$,
$61=4^2+5\cdot (-3)^2$, $521=21^2+5\cdot 4^2=(-11)^2+20^2$ and
$89=(-3)^2+5\cdot 4^2=5^2+8^2$. 

If $b=3$ then $d=13$ and we take $p=113,17,53,181,29,433,233$. We have
$113=10^2+13\cdot 1^2=(-7)^2+8^2$, $17=2^2+13\cdot 1^2=1^2+4^2$,
$53=1^2+13\cdot 2^2=(-7)^2+2^2$, $181=8^2+13\cdot (-3)^2$,
$29=4^2+13\cdot 1^2$, $433=(-15)^2+13\cdot 4^2=17^2+(-12)^2$ and
$233=5^2+13\cdot 4^2=13^2+8^2$. 

If $b=4$ then $d=5$, same as for $b=1$, so we can take again
$p=281,41,29,109,61,521,89$. 

(II) $b=2,6,10,14$ cover all cases (1)-(4) of (II). For each of them
we list four primes $P$ that cover all cases a-d from 2.16(II).

If $b=2$ then $d=2$ and we take $p=73,17,113,41$. We have
$73=(-7)^2+2\cdot (-6)^2=(-3)^2+8^2$, $17=(-3)^2+2\cdot 2^2=1^2+4^2$,
$113=9^2+2\cdot 4^2=(-7)^2+8^2$ and $41=(-3)^2+2\cdot 4^2=5^2+4^2$. 

If $b=6$ then $d=10$ and we take $p=89,41,281,241$. We have
$89=(-7)^2+10\cdot 2^2=5^2+8^2$, $41=1^2+10\cdot 2^2=5^2+4^2$,
$281=(-11)^2+10\cdot 4^2=5^2+16^2$ and $241=9^2+10\cdot
4^2=(-15)^2+4^2$. 

If $b=10$ then $d=26$ and we take $p=113,1297,1889,641$. We have
$113=(-3)^2+26\cdot 2^2=(-7)^2+8^2$, $1297=(-19)^2+26\cdot
(-6)^2=1^2+(-36)^2$, $1889=(-15)^2+26\cdot 8^2=17^2+40^2$ and
$641=(-15)^2+26\cdot 4^2=25^2+4^2$. 

If $b=14$ then $d=50$ and we take $p=281,641,881,809$. We have
$281=(-3)^2+50\cdot 2^2=5^2+16^2$, $641=21^2+50\cdot 2^2=25^2+4^2$,
$881=9^2+50\cdot 4^2=25^2+16^2$ and $809=(-3)^2+50\cdot 4^2=5^2+28^2$. 

(III) $b=8,16$ cover the cases (1) and (2) of (III). For each of them
we list six primes $p$ covering the cases a-f of 2.16(III).

If $b=8$ then $d=17$ and we take $p=53,141,593,409,1361,281$. We have
$53=(-6)^2+17\cdot 1^2=(-7)^2+2^2$, $141=9^2+17\cdot 2^2$,
$593=(-24)^2+17\cdot 1^2=(-23)^2+8^2$, $409=16^2+17\cdot
(-3)^2=(-3)^2+20^2$, $1361=(-33)^2+17\cdot 4^2=(-31)^2+20^2$ and
$281=(-3)^2+17\cdot 4^2=5^2+16^2$. 

If $b=16$ then $d=65$ and we take $p=101,269,1361,601,5009,1049$. We
have $101=(-6)^2+65\cdot 1^2$, $269=(-3)^2+65\cdot 2^2=13^2+10^2$,
$1361=36^2+65\cdot 1^2=(-31)^2+20^2$, $601=4^2+65\cdot
(-3)^2=5^2+24^2$, $5009=(-63)^2+65\cdot 4^2=65^2+(-28)^2$ and
$1049=(-3)^2+65\cdot 4^2=5^2+32^2$. 

But these pairs $b,p$ lie in the range of numbers already verified by
Sun. This concludes our proof. 

\section{Related results}

Let now $d>1$ be a square-free integer and let $\varepsilon\in\QQ (\sqrt d)$ be an
algebraic integer. Our problem is to determine $\varepsilon^{\frac{p-1}4}\mod p$
for primes $p$ satisfying $\legendre dp=\legendre{-1}p=1$. The answer will be
given in terms of $x,y,u,v\in\ZZ$ satisfying $p=f(x,y)=u^2+v^2$, where
$f(x,y)=ax^2+bxy+cy^2$ is a primitive positive definite quadratic form
of discriminant $b^2-4ac=-d$ or $-4d$, corresponding to $-d\ev 1\pmod
4$ or $-d\ev 2,3\pmod 4$, respectively. For a prime $q$, including
$q=\j$, we denote by $f_q$ the localized of $f$ at $q$. 

We proceed similarly as in [B, $\S$3]. First note that, in order to be
eligible, $f$ should represent elements $\ev 1\pmod 4$. In $\QQ_2$ a
number $\ev 1\pmod 4$ is in the square class of $1$ or $5$. Thus $f_2$
should represent $1$ or $5$. Assume this happens. We show that the
quadratic form $F(x,y,u,v)=f(x,y)-(u^2+v^2)$ is isotropic. By the
Hasse-Minkowski theorem we have to prove the statement locally. Since
both $f(x,y)$ and $u^2+v^2$ are positive definite $F_\j$ is
isotropic. If $q>2$ is prime and $q\nmid d$ then $F_\q$ is unimodular
so isotropic. If $q>2$ and $q\mid d$ then $F_q$ is isotropic because
$\det F_q=d\neq 1$ in $\QQ_q^\ti/(\QQ_q^\ti )^2$. At $q=2$ we have
that $f_2(x,y)$ represents $1$ or $5$ and $u^2+v^2$ represents both
$1$ and $5$ so $F_2$ isotropic. 

Take $(x_1,y_1,u_1,v_1)\in \QQ^4\setminus\{ (0,0,0,0)\}$ such that
$F(x_1,y_1,u_1,v_1)=0$. Since both $f(x,y)$ and $u^2+v^2$ are
anisotropic we have $f(x_1,y_1)=u_1^2+v_1^2=:a'\neq 0$. By multiplying
$x_1,y_1,u_1,v_1$ with a proper rational number we may assume that
$x_1,y_1\in\ZZ$ and $(x_1,y_1)=1$. Now $a'\in\ZZ$ and $a'$ is
represented by $u^2+v^2$ over $\QQ$ so it will be represented also
over $\ZZ$. So may assume that $u_1,v_1\in\ZZ$. Now $a'$ is
represented by $f$ primitively so we can write $f(x,y)=g(x',y')$,
where the mapping $(x,y)\mapsto (x',y')$ belongs to $SL(2,\ZZ )$ and
$g$ has the form $g(x',y')=a'x'^2+bx'y'+cy'^2$. We have
$b'^2-4a'c'=b^2-4ac=-d$ or $-4d$. 

If $p>2$ is a prime, $p\nmid a'd$ with $\legendre dp=\legendre{-1}p=1$ such that
$p=f(x,y)=u^2+v^2$ then
$$a'p=a'g(x',y')=(a'x'+\frac{b'}2y')^2+\frac{4a'c'-b'^2}4y'^2.$$ 
But $4a'c'-b'^2=d$ or $4d$, according as $-d\ev 1\pmod 4$ or
$-d\ev 2,3\pmod 4$. So $a'p=X^2+dY^2$, where $X=a'x'+\frac{b'}2y'$ and
$Y=\frac 12y'$ or $y'$, respectively. We also have
$a'p=(u_1^2+v_1^2)(u^2+v^2)=U^2+V^2$, where $U=u_1u-v_1v$ and
$V=u_1v+v_1u$. Note that $X,Y$ are linear combinations of $x',y'$ and
so of $x,y$ and $U,V$ are linear combinations of $u,v$. 

We have $a'p=X^2+dY^2=U^2+V^2$, which resembles the relations
$p=x^2+dy^2=u^2+v^2$ from $\S$2. Therefore we repeat the reasoning
from $\S$2 with $p,x\pm y\sqrt di,u\pm vi$ replaced by $a'p,X\pm Y\sqrt
di,U\pm Vi$. We will consider the fields $F=\QQ (\sqrt d)$,
$E=F(\mu_4)=\QQ (\sqrt d,i)$ and $L=E(\sqrt[4]{A_1})$, where $A_1=(X-Y\sqrt
di)(U+Vi)\1$. The analogue of Lemma 2.1 will hold. We keep the
notation $\s,\tau$ from the proof of Lemma 2.1 and define $\chi
:Gal(L/F)\z\mu_4$ by $\s^k\tau_l\mapsto i^k$. 

Same as for Lemma 2.2(i), we have
$$\chi\left(\legendre{\varepsilon,L/F}\p\right)\ev\varepsilon^{-\frac{p-1}4}\text{ so
}\varepsilon^{\frac{p-1}4}\ev\eta
:=\prod_{\q\neq\p}\chi\left(\legendre{\varepsilon,L/F}\q\right)$$
modulo $\P$. The primes $\q\neq\p$ where
$\chi\left(\legendre{\varepsilon,L/F}\q\right)$ may be $\neq 1$ are $\j$ and the
non-archimedian primes containing $2a'\varepsilon$. The challenge is to
calculate the Artin symbol at these primes. In a similar way one shows
that $\(\varepsilon^{\frac{p-1}4}\ev\eta'
:=\prod_{\q\neq\p}\chi\left(\legendre{\(\varepsilon,L/F}\q\right)$. The value of
$\varepsilon^{\frac{p-1}4}\mod p$ can be found from $\eta,\eta'$ as follows.

\blm If $\alpha,\beta\in\{\pm 1\}$ then:

(i) If $\eta=\alpha$, $\eta'=\beta$ then $\varepsilon^{\frac{p-1}4}\ev\frac 12(\alpha +\beta )+\frac
12(\frac{YV}{XU}\beta -\frac{YV}{XU}\alpha )\sqrt d\pmod p$. 

(ii) If $\eta=\alpha$, $\eta'=\beta i$ then $\varepsilon^{\frac{p-1}4}\ev\frac 12(\alpha -\frac VU\beta
)+\frac 12(\frac YX\beta -\frac{YV}{XU}\alpha )\sqrt d\pmod p$. 

(iii) If $\eta=i\alpha$, $\eta'=\beta$ then $\varepsilon^{\frac{p-1}4}\ev\frac 12(\beta -\frac
VU\alpha )+\frac 12(\frac{YV}{XU}\beta -\frac YX\alpha )\sqrt d\pmod p$. 

(iv) If $\eta=i\alpha$, $\eta'=i\beta$ then $\varepsilon^{\frac{p-1}4}\ev\frac 12(-\frac VU\alpha
-\frac VU\beta )+\frac 12(\frac YX\beta -\frac YX\alpha )\sqrt d\pmod p$.
\elm
\pf Let $\varepsilon^{\frac{p-1}4}=A+B\sqrt d$. Then  $\(\varepsilon^{\frac{p-1}4}=A-B\sqrt d$. It
follows that $A+B\sqrt d\ev\eta\pmod\P$ and $A-B\sqrt d\ev\eta'\pmod\P$ and
so $A\ev\frac 12(\eta +\eta')\pmod\P$ and $A\ev\frac 12(\frac 1{\sqrt d}\eta -\frac
1{\sqrt d}\eta')\pmod\P$. In the four cases (i)-(iv) of our lemma we get
that $A$ is $\ev\mod\P$ to $\frac 12(\alpha +\beta )$, $\frac 12(\alpha +\beta i)$, $\frac
12(\alpha i+\beta )$ or $\frac 12(\alpha i+\beta i)$, respectively, while $B$ is
$\ev\mod\P$ to $\frac 12(\frac 1{\sqrt d}\alpha -\frac 1{\sqrt d}\beta )$, $\frac 12(\frac 1{\sqrt
  d}\alpha -\frac i{\sqrt d}\beta )$, $\frac 12(\frac i{\sqrt d}\alpha -\frac 1{\sqrt d}\beta )$ or $\frac
12(\frac i{\sqrt d}\alpha -\frac i{\sqrt d}\beta )$, respectively. 

But $X-Y\sqrt di,U-Vi\in\P$ so $X\ev Y\sqrt d i$ and $U\ev Vi\pmod\P$. It
follows that $\frac i{\sqrt d}\ev -\frac YX\pmod\P$ and $i\ev -\frac VU\pmod\P$,
which by multiplication yield to $\frac 1{\sqrt d}\ev
-\frac{YV}{XU}\pmod\P$. So $A$ is $\ev\mod\P$ to $\frac 12(\alpha +\beta )$, $\frac
12(\alpha -\frac VU\beta )$, $\frac 12(\beta -\frac VU\alpha )$ or $\frac 12(-\frac VU\alpha -\frac VU\beta
)$, respectively, and $B$ is $\ev\mod\P$ to $\frac
12(\frac{YV}{XU}\beta -\frac{YV}{XU}\alpha )$, $\frac 12(\frac YX\beta -\frac{YV}{XU}\alpha )$, $\frac
12(\frac{YV}{XU}\beta -\frac YX\alpha )$ or $\frac 12(\frac YX\beta
-\frac YX\alpha )$,
respectively. But both sides of these congruences are $p$-adic
integers so the congruences also hold $\mod p$. Hence the
conclusion. \qed

Note that in the case of Conjectures 1-4 only the cases (i) and (iv)
of Lemma 3.1 can happen. (See the proof of Lemma 2.4.) This happens
because $\varepsilon\(\varepsilon =-1$ so
$\(\varepsilon^{\frac{p-1}4}=(-1)^{\frac{p-1}4}\varepsilon^{-\frac{p-1}4}$, which implies that
$\eta'=(-1)^{\frac{p-1}4}\eta\1$. In general, if $\varepsilon\(\varepsilon =N$ then
$\eta'\ev N^{\frac{p-1}4}\eta\pmod\P$. If $\legendre Np=1$ then
$N^{\frac{p-1}4}\ev\pm 1\pmod\P$ and so $\eta'=\pm\eta\1$ so $\eta,\eta'$
are either both of the type $\pm 1$ or both of the type $\pm i$. If
$\legendre Np=-1$ then $\eta,\eta'$ are of different type so the case (ii)
or (iii) of Lemma 3.1 applies. 
\vskip 3mm

{\bf Acknowledgement} I want to thank porfessor Zhi-Hong Sun for
introducing me to his conjectures.

\section*{References}

\hskip 6mm [B] C.N. Beli, {\it Two conjectures by Zhi-Hong Sun}, Acta
Arith. 137 (2009), 99-131.

[S1] Z.H. Sun, {\it Values of Lucas sequences modulo primes}, Rocky
Mountai J. Math. 33 (2003), 1123-1145.

[S2] Z.H Sun, {\it Quartic, octic residues and Lucas sequences},
J. Number Theory 129 (2009), 499-550.
\vskip 3mm

Institute of Mathematics ``Simion Stoilow'' of the Romanian Academy

21 Calea Grivitei Street, 010702-Bucharest, Sector 1, Romania
\end{document}